\documentclass[reqno]{amsart}

\usepackage[title]{appendix}
\usepackage[french,english]{babel}
\usepackage[T1]{fontenc}
\usepackage{hyperref}
\usepackage{comment}
\usepackage{palatino}
\usepackage[autostyle]{csquotes}
\usepackage{amsmath}
\usepackage{amsfonts,amssymb,amsmath,latexsym,mathrsfs,bbm,stmaryrd}
\usepackage[all]{xy}
\usepackage[usenames]{color}
\usepackage{epsfig}
\usepackage{graphicx}
\usepackage{subcaption}
\usepackage{float,enumerate}
\usepackage{amsthm}
\usepackage{thmtools}
\usepackage{thm-restate}
\usepackage[mathcal]{euscript}
\usepackage[nocompress]{cite}

\theoremstyle{plain}
\newtheorem{theorem}{Theorem}[section]

\newtheorem{assumption}[theorem]{Assumption}

\newtheorem{lemma}[theorem]{Lemma}

\theoremstyle{definition}
\newtheorem{definition}[theorem]{Definition}

\newtheorem{remark}[theorem]{Remark}

\numberwithin{equation}{section}

\renewcommand{\Im}{\operatorname{Im}}
\newcommand{\un}{1\mkern -4mu{\rm l}}
\renewcommand{\P}{\mathbb{P}}
\renewcommand\Re{\operatorname{Re}}
\renewcommand\Im{\operatorname{Im}}

\def\E{{\mathbb E}}

\def\i{{\rm i}}

\def\eps{\varepsilon}
\def\<{{\langle}}
\def\>{{\rangle}}

\DeclareMathOperator{\Tr}{tr}

\DeclareMathOperator{\dist}{dist}

\DeclareMathOperator{\Sparse}{Sparse}
\DeclareMathOperator{\supp}{supp}
\DeclareMathOperator{\Comp}{Comp}
\DeclareMathOperator{\Incomp}{Incomp}
\DeclareMathOperator{\Span}{Span}

\begin{document}

\title{Spectrum of Laplacian matrices associated with large random elliptic matrices}

\author{
Sean O'Rourke, Zhi Yin, and Ping Zhong
}

\address{
\parbox{\linewidth}{Sean O'Rourke, 
Department of Mathematics, 
University of Colorado,\\ 
Campus Box 395,
Boulder, CO 80309-0395, USA. \\
\texttt{sean.d.orourke@colorado.edu}}
}

\address{
\parbox{\linewidth}{Zhi Yin,
School of Mathematics and Statistics, 
Central South University,\\
Changsha, Hunan 410083, China.\\
\texttt{hustyinzhi@163.com}}
}

\address{
	\parbox{\linewidth}{Ping Zhong,
	Department of Mathematics and Statistics,
	University of Wyoming,\\
	Laramie, WY 82070, USA.\\
	\texttt{pzhong@uwyo.edu}}
	}

\date{\today}
\maketitle

\begin{abstract}
A Laplacian matrix is a square matrix whose row sums are zero.  
We study the limiting eigenvalue distribution of a Laplacian matrix formed by taking a random elliptic matrix and subtracting the diagonal matrix containing its row sums. Under some mild assumptions, we show that the empirical spectral distribution of the Laplacian matrix converges to a deterministic probability distribution as the size of the matrix tends to infinity. The limiting measure can be interpreted as the Brown measure of the sum of an elliptic operator and a freely independent normal operator with a Gaussian distribution.
\end{abstract}


\tableofcontents


\section{Introduction and main results}
Let $X$ be an $n\times n$ matrix with complex entries. The empirical spectral measure (ESM) $\mu_{X}$ of $X$ is a random probability measure defined as
\begin{equation}\label{eqn:ESD-X}
	    \mu_{X}=\frac{1}{n}\sum_{k=1}^{n}\delta_{\lambda_k(X)},
\end{equation}
where $\lambda_1(X), \lambda_2(X), \ldots,\lambda_n(X) \in \mathbb{C}$ are the eigenvalues of $X$ (counted with algebraic multiplicity) and $\delta_x$ is a point mass at $x$. If $X$ is an $n \times n$ Wigner random matrix, then its eigenvalues are real and $\mu_{X/\sqrt{n}}$ converges almost surely to the semicircular law as $n\rightarrow\infty$ \cite{book2010RMT}. If $X$ is an i.i.d. random matrix whose entries are independent and identically distributed (i.i.d.) copies of a complex random variable with mean zero and variance one, then the circular law says that $\mu_{X/\sqrt{n}}$ converges almost surely to the uniform measure on the unit disk \cite{TaoVu2010}. The elliptic random matrix model is a natural interpolation between these two random matrix models, initially introduced by Girko \cite{MR0781018,MR0901506,MR1483014,MR0816278,MR1354817,MR1373142,MR2999219,MR3068415,MR2215672}. The limiting ESM of an elliptic random matrix model is the uniform probability measure supported in an ellipsoid (see \cite{Naumov13} and \cite{NguyenORourke2015}). 

The purpose of this article is to study the limiting ESM of a random matrix of the form
\begin{equation}
	\label{eqn:defn-Laplacian}
	L=X-D,
\end{equation}
where $X$ is an elliptic random matrix and $D$ is the diagonal matrix whose diagonal entries are the row sums of $X$: 
\begin{equation} \label{eq:def:D}
   D_{ii}=\sum_{k=1}^{n}X_{ik}. 
\end{equation}
The matrix $L$ is called the \emph{Laplacian matrix formed from $X$}. 
If the entries $X_{ij}$ are nonnegative, the matrix $X$ can be interpreted as the adjacency matrix of a weighted oriented graph, and $L$ is the associated combinatorial Laplacian matrix. The matrix $L$ can be regarded as the infinitesimal generator of the continuous-time random walk on that graph. The spectrum of $L$ may be used to study the dynamics of the random walk on the graph. 

The case when $X$ is Wigner random matrix has been studied by Bryc, Dembo, and Jiang \cite{BrycDemboJiang2006}.

\begin{theorem}[Bryc--Dembo--Jiang] \label{thm:wigner}
Let $\xi$ be a real-valued random variable with mean zero and unit variance.  Let $X = (X_{ij})$ be an $n \times n$ random symmetric matrix whose entries $X_{ij}$, $n \geq j \geq i \geq 1$ are i.i.d. copies of $\xi$.  Let $L$ and $D$ be defined as in \eqref{eqn:defn-Laplacian} and \eqref{eq:def:D}.  Then the ESM of $L/\sqrt{n}$ converges weakly almost surely as $n \to \infty$ to the free convolution of the semicircular and standard normal distributions.  
\end{theorem}

The case when $X$ is a non-Hermitian matrix with i.i.d. entries was studied by Bordenave, Caputo, and Chafa\"{\i} \cite{BordenaveCaputoChafai2014markov}.  

\begin{theorem}[Bordenave--Caputo--Chafa\"{\i}] \label{thm:iid}
Let $\xi$ be a complex-valued random variable with mean zero and unit variance.  Let $X = (X_{ij})$ be an $n \times n$ matrix whose entries are i.i.d. copies of $\xi$.  Let $L$ and $D$ be defined as in \eqref{eqn:defn-Laplacian} and \eqref{eq:def:D}.  Then the ESM of $L/\sqrt{n}$ converges weakly in probability as $n \to \infty$ to a deterministic probability measure, which is the Brown measure of the sum of a circular operator and a $\ast$-free normal operator with a Gaussian distribution.  
\end{theorem}

Since elliptic random matrices interpolate between Wigner and i.i.d. matrices, our main results interpolate between and generalize these two results.  

\subsection{Laplacian matrices formed from elliptic random matrices}
In the Wigner random matrix model and the i.i.d. random matrix model, the $(i,j)$-entry $X_{ij}$ of $X$ is independent from all other entries, except possibly the $(j,i)$-entry $X_{ji}$. Elliptic random matrices interpolate between these two models.

\begin{definition}[The elliptic random matrix model with {\bf C0} condition]\label{defn:C0} 
Let $(\xi_1, \xi_2)$ be a random vector in $\mathbb{C}^2$ with mean zero and each entry having variance one. 
For each integer $n \geq 1$, we define the $n \times n$ random matrix $Z = (Z_{ij})_{i,j=1}^n$. We say the random matrix $Z$ is an \emph{elliptic random matrix model satisfying condition {\bf C0} with atom variables $(\xi_1, \xi_2)$} if the following three conditions hold:\\
\begin{enumerate}
	\item[{\bf C0-a:}] $\{Z_{ii}: 1\leq i \leq n \} \cup \{(Z_{ij}, Z_{ji}): 1\leq i<j \leq n\}$ is a collection of independent random elements;\\
    \item[{\bf C0-b:}] $\{(Z_{ij}, Z_{ji}): 1\leq i <j \leq n \}$ is a collection of i.i.d. copies of $(\xi_1, \xi_2)$;\\
    \item[{\bf C0-c:}] $\{Z_{ii}, 1\leq i \leq n \}$ is a collection of i.i.d. random variables with mean zero and finite variance.
\end{enumerate}
\end{definition}

We will assume the atom variable satisfy the following. 
\begin{definition}[$(\mu, \gamma)$-family]\label{defn:atom-variable}
Given parameters $0 \leq \mu \leq 1$ and $-1 < \gamma < 1,$ we say that the complex 
	random variables $(\xi_1, \xi_2)$ belongs to the $(\mu, \gamma)$-family if the following holds:
	\begin{enumerate}[\rm (1)]
		\item $\xi_1$ and $\xi_2$ have mean zero and variance one;
		\item $\mathbb{E} [(\Re \xi_1)^2] = \mathbb{E} [(\Re \xi_2)^2] =\mu$ and $\mathbb{E} [(\Im \xi_1)^2] = \mathbb{E} [(\Im \xi_2)^2] =1-\mu;$
		\item $\mathbb{E} [\Re \xi_1 \cdot \Re \xi_2] = \mu \gamma$ and $\mathbb{E} [\Im \xi_1 \cdot \Im \xi_2] = -(1-\mu) \gamma;$
		\item $\mathbb{E} [\Re \xi_i \cdot \Im \xi_j] = 0$ for any $i,j =1,2.$
	\end{enumerate}
\end{definition}
Note that if $(\xi_1, \xi_2)$ belongs to the $(\mu, \gamma)$-family, then $\mathbb{E} [\vert \xi_i \vert^2] =1, i=1,2$ and $\mathbb{E} [\xi_1 \xi_2] =\gamma.$ The covariance matrices are ${\rm Cov} (\Re \xi_i, \Im \xi_i) = \begin{bmatrix} \mu &
0 \\  0 & 1-\mu \end{bmatrix},$ $i=1,2.$ 

 When $Z$ is a random matrix that satisfies condition {\bf C0} with jointly Gaussian atom variables $(\xi_1, \xi_2)$ with $\mathbb{E} [\xi_1\xi_2] = \gamma$, the joint eigenvalue density can be calculated explicitly and the limiting ESM can be derived directly (see \cite{Gernotbook} and references therein). The limiting ESM is the uniform measure supported in the ellipsoid defined by 
 \begin{equation}\label{eqn:ellipsoid}
 \mathcal{E}_\gamma : = \left\{ (x, y) \in \mathbb{R}^2 : \frac{x^2}{(1+\gamma)^2} + \frac{y^2}{(1-\gamma)^2} \leq 1\right\}.
 \end{equation}
 In \cite{NguyenORourke2015}, Nguyen and the first author proved that, for any elliptic random matrix model satisfying condition {\bf C0} with atom variables from the $(\mu, \gamma)$-family, the 
 ESM of $Z/\sqrt{n}$ converges in probability as $n \rightarrow \infty$ to the uniform probability measure on the ellipsoid $\mathcal{E}_\gamma$. The real case was obtained by Naumov \cite{Naumov13} by assuming the atom variables have finite fourth moments. This universality result is the elliptic counterpart of the circular law established by Tao and Vu \cite{TaoVu2010} after a series of partial results; see \cite{Bordenave-Chafai-circular} for a detailed historic discussion. Additional spectral properties of elliptic random matrices have been studied in \cite{MR3403838,MR3357969,MR4565612,MR3540493,MR3744883,MR4388923,2013Low,MR4489824,MR4201596} and references therein.  
 
We now consider an elliptic random matrix $X$ satisfying only Conditions {\bf C0-a} and {\bf C0-b}. More precisely, given a random vector $(\xi_1, \xi_2)$ in $\mathbb{C}^2$, we assume $X = (X_{ij})$ is an $n \times n$ random matrix whose entries satisfy 
\begin{enumerate}[\rm (1)]
\item $\{X_{ii}: 1\leq i \leq n \} \cup \{(X_{ij}, X_{ji}): 1\leq i<j \leq n\}$ is a collection of independent random elements;
\item $\{(X_{ij}, X_{ji}): 1\leq i <j \leq n \}$ is a collection of i.i.d. copies of $(\xi_1, \xi_2)$;
\item the atom variables $(\xi_1, \xi_2)$ belongs to the $(\mu, \gamma)$-family,
\end{enumerate}
and there are no additional assumptions on the diagonal entries.  We use the notation $X$ instead of $Z$ for this matrix to distinguish that we are not assuming Condition {\bf C0-c}.  
We consider the Laplacian matrix $L=X-D$, where $D$ is the diagonal matrix obtained from the row sums of $X$ as defined in \eqref{eq:def:D}.  
We denote the rescaled matrix
\begin{equation}\label{eqn:def-M}
M = \frac{1}{\sqrt{n}}L=\frac{X-D}{\sqrt{n}}.
\end{equation}
In this paper, we show the following:
\begin{itemize}
\item the convergence of the ESM of $M$ as $n \to \infty$ to a deterministic limit, and 
\item a characterization of the limiting distribution in terms of free probability.
\end{itemize}



\subsection{Brown measure and elliptic operators}\label{subsec:brown}
We briefly review some basic concepts of free probability and refer the reader to \cite{MingoSpeicherBook} and the references therein for further details. 

Recall that a noncommutative probability space is a pair $(\mathcal{A}, \phi)$ where $\mathcal{A}$ is a von Neumann algebra and $\phi$ is a normal, faithful, tracial state on $\mathcal{A}$. 
Let $c$ be a Voiculescu's circular operator in $\mathcal{A}$ with variance one. Given $\gamma\in [-1,1]$, let $t_1=(1+ \gamma)/2, t_2=(1- \gamma)/2$ and let $s_{t_1}, s_{t_2}$ be two freely independent semicircular operators with variances $t_1, t_2$ respectively. 
The \emph{elliptic operator} $g_{\gamma}$ with parameter $\gamma$ is defined as
\[
   g_{\gamma}=s_{t_1}+\sqrt{-1}s_{t_2}, 
\]
where $\sqrt{-1}$ denotes the imaginary unit.  
It is an interpolation between the circular operator $c$ and the semicircular operator $s$. The only nonzero free cumulants of $g_{\gamma}$ are given by
\begin{equation}\label{eqn:free-cumulnts-y}
\kappa(y, y^*)=\kappa(y^*,y)=1, \; \kappa(y, y)=\kappa(y^*, y^*)=\gamma,
\end{equation}
where $y=g_{\gamma}$. The operator $g_{\gamma}$ includes the following operators as special cases:
\begin{itemize}
	\item if $\gamma=0$, then $y$ is a circular operator $c$ with variance one; 
	\item if $\gamma=1$, then $y$ is a semicircular operator $s$ with variance one. 
\end{itemize}
The elliptic operator can be twisted with a phase factor $e^{\sqrt{-1}\theta}$, $\theta \in [0, 2\pi]$ where the parameter $\gamma$ can be any complex number $|\gamma|\leq 1$, but we restrict to real parameter in this article. 

We use the standard notation $|x|=\sqrt{x^*x}$ to denote the square root of the positive operator $x^*x$ by functional calculus. 
We denote by $\widetilde{\mathcal{A}}$ the collection of operators affiliated with $\mathcal{A},$ and $\log^+(\mathcal{A})$ the collection of those $x\in\widetilde{A}$ such that
\[
\phi( \log^+(|x|)<\infty.
\]
For any $x\in \log^+(\mathcal{A})$, the Brown measure $\mu_x$ of $x$ is the probability measure defined as
\begin{equation}
\mu_x=\frac{1}{2\pi}\triangle_z \left( \phi( \log|x- z|) \right), \; z \in\mathbb{C},
\end{equation}
where $\triangle_z=\partial^2_{\Re z}+\partial^2_{\Im z}$ is taken in the distribution sense. 
When $x$ is a self-adjoint operator or a normal operator in the sense that $xx^*=x^*x$, then $\mu_x$ is the same as the spectral measure of $x$ with respect to $\phi$. 

\subsection{Main results}
Given an $n \times n$ matrix $A$, recall that $\lambda_1(A),\lambda_2(A), \ldots, \lambda_n(A) \in \mathbb{C}$ are the eigenvalues of $A$.  The singular values of $A$ are the eigenvalues of $\sqrt{A^\ast A}$, where $A^\ast$ is the conjugate transpose of $A$.  Let $s_1(A) \geq  s_2(A) \geq \cdots \geq s_n (A)$ denote the ordered singular values of $A$. The ESM of $A$ is defined as in \eqref{eqn:ESD-X}: 
\[
\mu_A : = \frac{1}{n} \sum_{k=1}^n \delta_{\lambda_k(A)}.
\]
We define the empirical measure $\nu_A$ constructed from the singular values of $A$ by 
\begin{equation}
\nu_A : = \frac{1}{n} \sum_{k=1}^n \delta_{s_k(A)}.
\end{equation}
Denote $\mathbb{C}_+ : = \{z \in \mathbb{C}: {\Im} (z)>0 \}.$ For any probability measure $\nu$ on $\mathbb{R},$ its Stieltjes transform $S_{\nu}: \mathbb{C}_+ \rightarrow \mathbb{C}_+$ is given by 
\begin{equation}
S_\nu (z) = \int \frac{1}{t-z} d\nu(t).
\end{equation}
We denote by $\widetilde{\nu}$ the symmetrization of $\nu$, defined for any Borel set $\mathcal{B} \subseteq \mathbb{R}$ by
\begin{equation}
\widetilde{\nu}(\mathcal{B}) = \frac{\nu(\mathcal{B}) + \nu(-\mathcal{B})}{2}, 
\end{equation}
where $-\mathcal{B} = \{ -b : b \in \mathcal{B}\}$.  

We begin with the limit law for the singular values of $M - z$ for any $z \in \mathbb{C}$. In what follows, $G$ is a Gaussian random variable on $\mathbb{C}$ with law $\mathcal{N}(0, K_\mu)$. The real and imaginary parts of $G$ are centered real Gaussian random variables with covariance matrix $K_\mu =\begin{bmatrix} \mu &
0 \\  0 & 1-\mu \end{bmatrix}.$

\begin{theorem}\label{thm:singular-value}
Let $X$ be an $n \times n$ elliptic random matrix satisfying Condition {\bf C0-a} and {\bf C0-b} with atom variables belonging to the $(\mu, \gamma)$-family. Set $M=(X-D)/\sqrt{n}$ as defined in \eqref{eqn:def-M}. 
For every $z \in \mathbb{C},$ there exists a deterministic probability measure $\nu_z$ on $[0,\infty)$ such that in probability,
 $$\nu_{M -z} \rightarrow \nu_z \; \text{weakly as} \; n \rightarrow \infty.$$
Moreover, the limit law $\nu_z$ is characterized as
\begin{equation}
  \label{eqn:nu-1.12}
      \nu_z=\mu_{|a+g_\gamma-z|},
\end{equation}
where $g_\gamma$ is an elliptic operator with parameter $\gamma$ and $a$ is a Gaussian distributed normal operator with law $\mathcal{N}(0, K_\mu)$, freely independent from $g_\gamma$. 
\end{theorem}

\begin{remark}
We remark that for $X$ being a random matrix with i.i.d. entries, the convergence of $\nu_{M-z}$ holds with probability $1;$ see \cite[Theorem 1.1]{BordenaveCaputoChafai2014markov}. Since the matrix $M-z$ has independent rows, by the concentration of measure phenomenon for matrices with independent rows, one can
show that $\nu_{M-z} - \mathbb{E} \nu_{M-z}$ converges weakly almost surely to zero. However, for the elliptic case, our techniques can only establish convergence in probability; see Lemma \ref{lem:epsilon5} below for details. 
\end{remark}

The next result concerns the eigenvalues of $M.$ We will make the following further assumption for the atom variables $(\xi_1, \xi_2).$
	\begin{assumption} \label{assump:atom}
	We say the atom variables $(\xi_1, \xi_2)$ satisfy Assumption \ref{assump:atom} if there exist $\eps_0, \delta_0 > 0$ so that
	almost surely 
	\begin{equation} \label{assump:cond}
		\sup_{z \in \mathbb{C}} \P (|\xi_1 - z| \leq \delta_0 | \xi_2) \leq 1- \eps_0 \; \text{\rm and} \;
		\sup_{z \in \mathbb{C}} \P (|\xi_2 - z| \leq \delta_0 | \xi_1) \leq 1- \eps_0. 
	\end{equation}
	\end{assumption}
	
Here, the conditional probabilities in \eqref{assump:cond} should be interpreted as regular conditional probabilities (see, for example, Section 4.1.3 in \cite{MR3930614} or \cite{MR1484954}). 
Condition \eqref{assump:cond} implies that knowing the value of $\xi_2$ does not almost surely determine $\xi_1$ (and vice versa).  In particular, this condition implies that
	\begin{equation} \label{eq:assump:noncond}
		\max_{i = 1, 2} \sup_{z \in \mathbb{C}} \P (|\xi_i - z| \leq \delta_0) \leq 1- \eps_0. 
	\end{equation}
	
There are many examples of atom variables $(\xi_1, \xi_2)$ satisfying Assumption \ref{assump:atom}.  For example, if $\xi_1$ and $\xi_2$ are real-valued jointly Gaussian random variables with mean zero, unit variance, and correlation $\gamma = \E[\xi_1 \xi_2]$ satisfying $|\gamma| < 1$, then $(\xi_1, \xi_2)$ satisfy Assumption \ref{assump:atom}.  

\begin{theorem}\label{thm:eigenvalue}
Let $X$ be an $n \times n$ elliptic random matrix satisfying Conditions {\bf C0-a} and {\bf C0-b} with atom variables belonging to the $(\mu, \gamma)$-family and satisfying Assumption \ref{assump:atom}. Set $M=(X-D)/\sqrt{n}$ as defined in \eqref{eqn:def-M}. 
Let $\mu$ be the probability measure on $\mathbb{C}$ defined by 
\begin{equation}
\mu = \frac{1}{2\pi} \triangle_z U \; \text{with} \; U(z)=  \int_0^\infty \log (t) d \nu_z(t),
\end{equation}
where the Laplacian $\triangle_z = \partial_{\Re z}^2+ \partial_{\Im z}^2$ is taken in the distribution sense, and $\nu_z$ is as in Theorem \ref{thm:singular-value}. Then, in probability, 
\begin{equation*}
\mu_M \rightarrow \mu \; \text{weakly as} \; n \rightarrow \infty.
\end{equation*}
In other words, the limiting distribution of $\mu_M$ is the Brown measure $\mu=\mu_{a+g_{\gamma}}$. 
\end{theorem}

\begin{remark}
Theorems \ref{thm:singular-value} and \ref{thm:eigenvalue} can be extended to the case when the entries of $X$ do not have mean zero; see Subsection \ref{sub:non-centered} for details.
\end{remark}

Our results show that free probability theory provides a natural candidate for the limiting distribution. If $X$ is an elliptic random matrix with jointly Gaussian atom variables $(\xi_1, \xi_2)$ and $\widetilde{D}$ is a diagonal matrix independent of $X$, then by the asymptotic freeness result of Voiculescu \cite{Voiculescu1991}, the sequence of random matrices $ \widetilde{D} + X/\sqrt{n}$ converges in $*$-moments to $x+g_\gamma$ in $(\mathcal{A},\phi)$, where $\widetilde{D}$ converges to $x$ in $*$-moments and $x$ is freely independent from $g_\gamma$. Theorem \ref{thm:eigenvalue} implies a similar result; however, in our setting, the problem is more difficult because $D$ depends heavily on $X$. Our results also interpolate between and generalize Theorems \ref{thm:wigner} and \ref{thm:iid}.  

 We now turn to some further properties of $\mu_{a+ g_\gamma}$. We define the following open set 
\begin{equation}
\label{defn:Xi-t}
 \Xi =\left\{ \lambda \in \mathbb{C} : \mathbb{E} \left[\vert G- \lambda\vert^{-2}\right]>1 \right\}.
\end{equation}
For $\lambda \in \Xi$, let $w(\lambda)$ be a function of  $\lambda$ taking positive values such that
\[
  \mathbb{E} \left[ ( \vert G -\lambda\vert^2+w(\lambda)^2)^{-1} \right]=1
\]
and $w(\lambda)=0$ for $\lambda \in\mathbb{C}\backslash\Xi$. 
We denote
\[
  \Phi_{\gamma} (\lambda)= \lambda + \gamma \cdot \mathbb{E} \left[ (\bar{\lambda} - \overline{G}) \left( \vert \lambda- G \vert^2+w(\lambda)^2 \right)^{-1}\right]. 
\]
Let $c$ be a circular operator freely independent from $a$ in $(\mathcal{A},\phi)$. 
It is known \cite{BordenaveCaputoChafai2014markov} that ${\rm supp} (\mu_{a+c}) = \left\{ \lambda \in \mathbb{C} : \mathbb{E} \left[\vert G- \lambda\vert^{-2}\right] \geq 1 \right\}$ (see also \cite[Section 4]{Zhong2021Brown_ctgt}). Moreover, the density function of $\mu_{a+c}$ has an explicit formula which is strictly positive in the set $\Xi$. As a corollary of the results in \cite{Zhong2021Brown_ctgt, BelinschiYinZhong2021Brown}, we have the following description of the measure $\mu_{a+g_\gamma}.$

\begin{theorem}[Section 5 of \cite{Zhong2021Brown_ctgt}] \label{thm:prop-Brown-sum}
 The Brown measure of $a+g_\gamma$ is the push-forward measure of the Brown measure of $a+c$ by the map $\lambda \mapsto \Phi_{\gamma}(\lambda)$. That is, if $\mathcal{B}$ is an arbitrary Borel measurable set in $\mathbb{C}$, then 
	\begin{equation*}
	  \mu_{a+g_{\gamma}}(\mathcal{B})=\mu_{a+c}(\Phi_{\gamma}^{-1}(\mathcal{B})).
	\end{equation*}
 \end{theorem}

Note that $\phi [f(a)] = \mathbb{E} [f(G)]$ for any measurable function $f$. If $K_\mu$ is degenerate ($\mu=0$ or $1$) and $\supp(G)\subset\mathbb{R}$, then the boundary of $\supp(\mu_{a+c})$ is the set
\[
\left\{ z=z_1+\sqrt{-1}z_2 : z_1, z_2 \in \mathbb{R}, \int_\mathbb{R} \frac{1}{(u-z_1)^2+z_2^2}d\mu_G(u) =1\right\},
\]
where $\mu_G$ is the standard real Gaussian distribution. In this case, the support of $\mu_{a+c}$ is symmetric with respect to the $x$-axis, and the boundary set in the upper-half plane can be identified as the graph of a function on $\mathbb{R}$. See \cite[Section 6.1]{Zhong2021Brown_ctgt} for details. Moreover, the pushforward map $\lambda\mapsto \Phi_{\gamma}(\lambda)$ is a homeomorphism of $\mathbb{C}$ by \cite[Section 6.2]{Zhong2021Brown_ctgt}. If $K_\mu$ is invertible, then $\Xi=\mathbb{C}$ and $\supp(\mu_{a+c})=\supp(\mu_{a+g_\gamma})=\mathbb{C}$. 

We end this subsection with the following remark:
\begin{remark}
As mentioned before, one of our motivations is that our results interpolate between and generalize the results for Wigner and i.i.d. random matrices. Importantly, free probability 
theory provides a unifying framework for studying the limiting ESM of those random matrix models. 

Another motivation for our work comes from studying the limiting ESM of sparse random matrix models. Sparse matrices appear in statistics, computer science, wireless communications, quantum physics, and many other fields. In contrast to the non-sparse case, the rigorous study of sparse matrix models, especially in the non-Hermitian case, emerged much later. For instance, the sparse circular law was proved by Basak and Rudelson in 2019 \cite{BR2019}. For sparse elliptic random matrices, the limiting ESM was studied in some physical literature \cite{NM-PRL, MNR-JPA, Baron2023}. However, the convergence of the ESM still lacks rigorous proof. Although we are dealing with the non-sparse model, some of our techniques are applicable to the sparse case.
\end{remark}


\subsection{Outline of the proof}

Following the standard Hermitian reduction trick \cite{Bordenave-Chafai-circular} for non-Hermitian random matrices, the proofs for Theorem \ref{thm:singular-value} and Theorem \ref{thm:eigenvalue} have two main parts:
\begin{enumerate}
	\item identify a family of nonrandom probability measures $\{\nu_z\}_{z\in\mathbb{C}}$, such that for almost all $z\in\mathbb{C}$, $\nu_{M-z}\rightarrow \nu_z$ weakly in probability as $n\rightarrow\infty$;
	\item show that $\log (\cdot)$ is uniformly integrable for $\{\nu_{M -z}\}_{n\geq 1}$.
\end{enumerate}
For part (1), we will consider the following $2n \times 2n$ Hermitian matrix:
\begin{equation*}
H(z) := \begin{bmatrix} 0 &
M-z \\  (M-z)^* & 0 \end{bmatrix}.
\end{equation*} 
Denote by $s_i(M-z)$ the singular values of $M-z$.  
The eigenvalues of $H(z)$ are exactly the pairs of singular values $\pm s_1(M-z), \ldots, \pm s_n(M-z)$. Hence, for $\eta \in \mathbb{C}_+$, the Stieltjes transform of $\widetilde{\nu}_{M-z}$ is given by
\begin{equation}
\label{eqn:intro-resolvent-main}
S_{\widetilde{\nu}_{M-z}} (\eta)= \frac{1}{2n} {\rm Tr} \left[ (H(z)- \eta)^{-1} \right].
\end{equation}
We derive some approximate fixed point equation satisfied by the limit of the $2\times 2$ block-wise Stieltjes transform of $H(z)$. Since the entries of $M$ are not i.i.d., 
we have to substantially modify the strategy of proof for the i.i.d. case in 
\cite{BordenaveCaputoChafai2014markov}. 
The estimate of the resolvent $(H(z)- \eta)^{-1}$ utilizes some techniques developed in \cite{2013Low, NguyenORourke2015} for the elliptic model and \cite{Tikhom2023} for the Hermitian model.

For part (2), we study a lower bound for the least singular value of $M-z$. The main challenge in establishing this lower bound is that the entries of $M$ are highly dependent on one another due to the diagonal entries containing the row sums.  After this bound is established, uniform integrability will follow from standard arguments. We adapt the argument of \cite[Proposition 3.4]{BordenaveCaputoChafai2014markov} and also a general approach for Hermitian matrices \cite{RV2008adv}, but some extra work is needed to deal with the dependence amongst the entries.

The rest of the paper is organized as follows. Section \ref{section:pre} collects some preliminary results. 
Section \ref{sec:approx} derives an approximate fixed point equation for the diagonal $2 \times 2$ block of the resolvent and gives the proof of Theorem \ref{thm:singular-value}. Section \ref{sec:lsv} is devoted to the bound on the least singular value, and Section \ref{sec:mslv} is for the moderately small singular values. In the last section, we finish the proof of Theorem \ref{thm:eigenvalue}. 

\subsection{Notation}
For $x=x(n)$ and $y=y(n)$, we use notation $|x|=O(y)$ to denote the bound $x \leq C y$ for all sufficiently large $n$ and for some universal constant $C$. Notation $x= o(y)$ means that $x/y \rightarrow 0$ as $n \rightarrow \infty.$ 

We write a.a. for Lebesgue almost all. We use $\sqrt{-1}$ to denote the imaginary unit and reserve $i$ as an index. 
 
 We denote an elliptic operator with parameter $\gamma$ by $g_{\gamma}$ and denote $G$
 by a Gaussian random variable on $\mathbb{C}$ with law $\mathcal{N}(0,K_\mu)$. 
 
 	For a finite set $S$, $|S|$ denotes its cardinality.  We use $[n]$ to denote the discrete interval $\{1, \ldots, n\}$.  For two random variables $X$ and $Y$, we write $X \stackrel{d}{=} Y$ to denote equality in distribution.  
	
	We let $\mathcal{M}_n(R)$ be the set of $n \times n$ matrices with entries in a ring $R$.  For any $n \times n$ matrix $A \in \mathcal{M}_n(\mathbb{C})$, we let $\lambda_1(A), \ldots, \lambda_n(A) \in \mathbb{C}$ denote the eigenvalues of $A$ (counted with algebraic multiplicity).  The ordered singular values of $A$ are denoted $s_1(A) \geq \cdots \geq s_n(A)$.  We let $\un_n$ denote the $n \times n$ identity matrix. Sometimes when its clear, e.g., $M - z \un_n$, we do not write the identity matrix, e.g.,  $M - z$.  
 		
	For a vector $v$, $\|v\|$ denotes its Euclidean norm.  For a matrix $A$, $A^{\mathrm{T}}$ is its transpose, $A^{\ast}$ is its conjugate transpose, and $\|A\|$ is its spectral norm (i.e., $\|A\|$ is the largest singular value of $A$).  
 		
 	For a vector $v = (v_i)_{i=1}^n \in \mathbb{C}^n$ and a subset $I \subset [n]$, we let $v_I = (v_i)_{i \in I} \in \mathbb{C}^I$.  Similarly, for an $m \times n$ matrix $A = (A_{ij})_{i \in [m], j \in [n]}$ and $I \subset [m], J \subset [n]$, we define $A_{I \times J} = (A_{ij})_{i \in I, j \in J}$.  
 		
 	For a vector $v \in \mathbb{C}^n$ and a subspace $H$ of $\mathbb{C}^n$, we write $\dist(v, H)$ to denote the Euclidean distance from $v$ to $H$, i.e., $\dist(v, H) = \inf_{u \in H} \|v - u\|$.


\section{Truncation and centralization}\label{section:pre}
We begin the proof of our main results by truncating the entries of the random matrix $M$.  In this section, we also show how our main results can be extended to the case where the entries have non-zero mean.  
\subsection{Truncation of the atom variables}

We define 
\begin{equation*}
\xi'_i = \xi_i \un_{\{\vert \xi_i \vert \leq \varepsilon \sqrt{n} \}}, \qquad i =1, 2, 
\end{equation*}
where $\varepsilon=\varepsilon(n) \rightarrow 0$ such that 
\begin{equation}\label{eq:truncation-1}
\lim_{n \rightarrow \infty} \mathbb{E} \left[ \vert \xi_i \vert^2 \un_{\{\vert \xi_i \vert^2 \geq \varepsilon^2 n \}} \right] =0.
\end{equation}
We also denote, for $i=1,2$, 
\[
 \hat{\xi}_i = \xi_i'- \mathbb{E} [\xi_i'], \qquad \hat{\gamma} = \mathbb{E} \left[ \hat{\xi}_1 \hat{\xi}_2 \right]
\]
and
\begin{equation*}
   \xi_i''=\xi_i-\xi_i'=\xi_i \un_{\{\vert \xi_i \vert > \varepsilon \sqrt{n} \}}.
\end{equation*}

\begin{lemma}[Truncation of the atom variables]\label{lem:truncation}
Let $\hat{X}$ be a $n \times n$ random matrix that satisfies Conditions {\bf C0-a} and {\bf C0-b} with atom variables $(\hat{\xi}_1, \hat{\xi}_2),$ and let $\hat{X}_{ii} = X_{ii}, i=1, \ldots, n.$ 
Let $\hat{D}$ be a diagonal matrix with $\hat{D}_{ii} = \sum_{k=1}^n \hat{X}_{ik}, i=1, \ldots, n.$ Set $\hat{M} = (\hat{X} - \hat{D})/\sqrt{n},$ 
then we have 
\begin{equation}
\mathbb{E} \nu_{\hat{M}-z} - \mathbb{E} \nu_{M-z} \rightarrow 0 \; \text{weakly and} \; \vert \gamma - \hat{\gamma} \vert \rightarrow 0
\end{equation}
as $n \rightarrow \infty.$ 
\end{lemma}

\begin{proof}
From the Hoffman-Wielandt inequality \cite{HWineq} 
\begin{equation*}
\begin{split}
\frac{1}{n}\mathbb{E}  \sum_{k=1}^n \vert s_k (\hat{M}-z) & -s_k (M-z) \vert^2  \leq \frac{1}{n} \mathbb{E}  \sum_{i, j=1}^n \vert \hat{M}_{ij} - M_{ij} \vert^2\\
& \leq \frac{2}{n^2}\mathbb{E} \left( \sum_{i, j=1}^n \vert \hat{X}_{ij} - X_{ij} \vert^2 + \sum_{i=1}^n \vert \hat{D}_{ii} - D_{ii} \vert^2 \right).
\end{split}
\end{equation*}
We have
\begin{equation}\label{eqn:2.7-in-proof}
\begin{split}
\frac{1}{n^2}\mathbb{E} \left( \sum_{i, j=1}^n \vert \hat{X}_{ij} - X_{ij} \vert^2 \right) & = \frac{1}{n^2}\mathbb{E} \left( \sum_{i\neq j} \vert \hat{X}_{ij} - X_{ij} \vert^2 \right) \\
&\leq \frac{1}{n^2} \sum_{i<j}\mathbb{E}|\hat{\xi}_1-\xi_1|^2 +\frac{1}{n^2} \sum_{i>j}\mathbb{E}|\hat{\xi}_2-\xi_2|^2.
\end{split}
\end{equation}
Note that, for $i=1,2$, 
\[
   \mathbb{E}(\xi_i')+\mathbb{E}(\xi_i'')=\mathbb{E}(\xi_i)=0.
\]
Hence, we can rewrite
\[
  \mathbb{E}|\hat{\xi}_i-\xi_i|^2=\mathbb{E}| \xi_i'-\xi_i-\mathbb{E}(\xi_i')  |^2
   =\mathbb{E}|\xi_i''-\mathbb{E}(\xi_i'')|^2
    \leq 2\mathbb{E}|\xi_i''|^2.
\]
We can continue to estimate
\begin{equation*}
  \begin{split}
    \text{Right-hand side of \eqref{eqn:2.7-in-proof}}
         &\leq \frac{1}{n^2}\sum_{i<j}\mathbb{E}|\xi_1''-\mathbb{E}(\xi_1'')|^2
           +\frac{1}{n^2}\sum_{i>j}\mathbb{E}|\xi_2''-\mathbb{E}(\xi_2'')|^2\\
        &\leq \frac{2}{n^2}\sum_{i<j}\mathbb{E}|\xi_1''|^2
          +\frac{2}{n^2}\sum_{i>j}\mathbb{E}|\xi_2''|^2
  \end{split}
\end{equation*}
which tends to $0$ as $n \rightarrow \infty$ by by \eqref{eq:truncation-1}.

Since the entries in the $i$-th row $\{X_{ik}\}_{i=1}^n$ are independent, by a similar calculation, we have 
\begin{equation*}
\begin{split}
\frac{1}{n^2}\mathbb{E}& \left( \sum_{i=1}^n \vert \hat{D}_{ii} - D_{ii} \vert^2 \right)  = \frac{1}{n^2}\mathbb{E} \left( \sum_{i=1}^n \vert \sum_{k\neq i}  \hat{X}_{ik}- \sum_{k\neq i} X_{ik} \vert^2 \right) \\
&\leq \frac{1}{n^2}\sum_{i=1}^n\sum_{k\neq i}\mathbb{E}\vert \hat{X}_{ik}- X_{ik}\vert^2 \\
&\leq \frac{1}{n^2} \sum_{i<j}\mathbb{E}|\hat{\xi}_1-\xi_1|^2 +\frac{1}{n^2} \sum_{i>j}\mathbb{E}|\hat{\xi}_2-\xi_2|^2\\
&\leq \frac{2}{n^2}\sum_{i<j}\mathbb{E}|\xi_1''|^2
          +\frac{2}{n^2}\sum_{i>j}\mathbb{E}|\xi_2''|^2
  \end{split}
\end{equation*}
which tends to $0$ as $n \rightarrow \infty$ by \eqref{eq:truncation-1}.
Hence, we deduce that 
\begin{equation}\label{eqn:2.8-in-proof}
\frac{1}{n}\mathbb{E}  \sum_{k=1}^n \vert s_k (\hat{M}-z) -s_k (M-z) \vert^2 \rightarrow 0.
\end{equation}
Since the left-hand side of \eqref{eqn:2.8-in-proof} bounds the square of the expected Wasserstein $W_2$ coupling distance between $\nu_{\hat{M}-z}$ and $\nu_{M-z}$, this yields that
\[
\mathbb{E} \nu_{\hat{M}-z} - \mathbb{E} \nu_{M-z} \rightarrow 0
\]
as $n \rightarrow \infty$.

We note that \[
  \gamma=\mathbb{E}(\xi_1\xi_2)=\mathbb{E}( (\xi_1'+\xi_1'')(\xi_2'+\xi_2''))
\]
and
\[
\hat{\gamma}=\mathbb{E}(\hat{\xi}_1\hat{\xi}_2)=
    \mathbb{E}(\xi_1'\xi_2')-\mathbb{E}(\xi_1')\mathbb{E}(\xi_2').
\]
Recall that $\xi_2''+ \xi_2'= \xi_2$. 
We obtain
\begin{equation} \label{eqn:2.9-in-proof}
  \begin{split}
    \gamma-\hat{\gamma}
      &= \mathbb{E}(\xi_1''\xi_2') + \mathbb{E}(\xi_1'' \xi_2'')+\mathbb{E}(\xi_1'\xi_2'')+\mathbb{E}(\xi_1')\mathbb{E}(\xi_2')\\
      &=\mathbb{E}(\xi_1''\xi_2)+\mathbb{E}(\xi_1'\xi_2'')+\mathbb{E}(\xi_1')\mathbb{E}(\xi_2')\\
      &=\mathbb{E}(\xi_1''\xi_2)+\mathbb{E}(\xi_1'\xi_2'')-\mathbb{E}(\xi_1')\mathbb{E}(\xi_2''),
  \end{split}
\end{equation}
where we used $\mathbb{E} (\xi_2') + \mathbb{E}(\xi_2'') =0$ for the second equality.
By the Cauchy-Schwarz inequality, we have
\[
  \vert\mathbb{E}(\xi_1''\xi_2)\vert^2\leq \mathbb{E}\vert\xi_1''\vert^2 \mathbb{E}\vert\xi_2\vert^2\leq \mathbb{E}\vert\xi_1''\vert^2
    =\mathbb{E}\vert\xi_1\vert^2 \un_{\{\vert \xi_1 \vert > \varepsilon \sqrt{n} \}}.
\]
We note that, for $i=1, 2$, 
\begin{equation*}
{\rm Var} (\hat{\xi}_i) \leq \mathbb{E} \left[ \vert \xi_i \vert^2 \un_{\{ \vert \xi_i \vert \leq \varepsilon \sqrt{n} \}} \right] \leq 1.
\end{equation*}
Hence, we have 
\[
  \begin{split}
    \vert \mathbb{E}(\xi_1'\xi_2'')-\mathbb{E}(\xi_1')\mathbb{E}(\xi_2'')\vert^2
      &=\vert \mathbb{E} ( (\xi_1'-\mathbb{E}(\xi_1')) (\xi_2''-\mathbb{E}(\xi_2'')))   \vert^2\\
      &\leq {\rm Var} (\hat{\xi}_1) \mathbb{E}\vert \xi_2''-\mathbb{E}(\xi_2'') \vert^2\\
      &\leq 2\mathbb{E}\vert \xi_2'' \vert^2
      =2\mathbb{E}\vert\xi_2\vert^2 \un_{\{\vert \xi_2 \vert > \varepsilon \sqrt{n} \}}.
  \end{split}
\]
Therefore, by \eqref{eqn:2.9-in-proof}, we have 
\begin{equation*}
  \begin{split}
    \vert \gamma-\hat{\gamma}\vert^2
   & \leq 2 \left( \vert\mathbb{E}(\xi_1''\xi_2)\vert^2+\vert \mathbb{E}(\xi_1'\xi_2'')-\mathbb{E}(\xi_1')\mathbb{E}(\xi_2'')\vert^2 \right)\\
    &\leq 2\mathbb{E}\vert\xi_1\vert^2 \un_{\{\vert \xi_1 \vert > \varepsilon \sqrt{n} \}}
      +4\mathbb{E}\vert\xi_2\vert^2 \un_{\{\vert \xi_2 \vert > \varepsilon \sqrt{n} \}}
  \end{split}
\end{equation*}
which tends to $0$ as $n\rightarrow\infty$ by by \eqref{eq:truncation-1}.
This finishes the proof. 
\end{proof}

In summary, for the remainder of the proof of Theorem \ref{thm:singular-value}, we will assume that 
\begin{enumerate}
\item $\mathbb{E} [\xi_i] =0,$ and $\mathbb{E} \vert \xi_i \vert^2 \leq 1, i=1,2;$
\item $\P (\max\{ |\xi_1|, |\xi_2| \} \geq \kappa (n)) =0,$
where $\kappa(n) = \eps\sqrt{n} = o(\sqrt{n}).$
\end{enumerate}

\subsection{Extension to non-centered entries}\label{sub:non-centered} 
In this section, we briefly note how our results can be extended to the case when the entries of the random matrix $X$ have non-zero mean.  Let $(\chi_1, \chi_2)$ be atom variables with expectation $m=\mathbb{E}(\chi_1)=\mathbb{E}(\chi_2)$ such that $(\chi_1-m, \chi_2-m)$ belongs to the $(\mu,\gamma)$-family defined in Definition \ref{defn:atom-variable}.  Suppose $X$ is a non-centered elliptic random matrix model satisfying the following conditions:
\begin{enumerate}[\rm (1)]
	\item $\{X_{ii}: 1\leq i \leq n \} \cup \{(X_{ij}, X_{ji}): 1\leq i<j \leq n\}$ is a collection of independent random elements,
	\item $\{(X_{ij}, X_{ji}): 1\leq i <j \leq n \}$ is a collection of i.i.d. copies of $(\chi_1, \chi_2),$
	\item $\{X_{ii}, 1\leq i \leq n \}$ is a collection of i.i.d. random variables. 
\end{enumerate}

Since we will be considering the matrix $L = X - D$, where $D$ is defined as in \eqref{eq:def:D}, without loss of generality, we assume the diagonal entries $X_{ii}$ of $X$ are zero.  
Consider the rescaled matrix
\[
   M=\frac{L+(n-1)m\un_n}{\sqrt{n}}=\frac{X}{\sqrt{n}}-\frac{D-((n-1)m)\un_n}{\sqrt{n}}.
\]
Note that $\mathbb{E}(M)=\frac{mJ}{\sqrt{n}}-\frac{m\un_n}{\sqrt{n}}$, where $J$ denotes the matrix with all entries equal to $1$. We set
\[
  \bar{M}=M-\mathbb{E}(M), \qquad  \bar{M}'=M-\frac{mJ}{\sqrt{n}}.
\]
Denote by $d_L (\cdot, \cdot)$ the L\'{e}vy distance between probability measures. By the norm inequality \cite[Theorem A. 45]{BaiBook}, we have 
\[
   d_L (\nu_{\bar{M}-z}, \nu_{\bar{M}'-z})\leq \Vert \bar{M}- \bar{M}' \Vert \leq \frac{|m|}{\sqrt{n}}.
\]
Moreover, by standard perturbation inequalities, we have
\[
   \left\vert \int f d\nu_{\bar{M}'-z}-\int f d\nu_{{M}-z} \right\vert \leq \Vert f\Vert_{\text{BV}}\frac{\text{rank}(M-\bar{M}')}{n}\leq \frac{\Vert f\Vert_{\text{BV}}}{n}.
\]
Hence, we can conclude that the eigenvalue distributions of $M$ and its centralization $\bar{M}$ both converge to the Brown measure of $a+g_\gamma$ under Assumption \ref{assump:atom}.



\section{Convergence of the singular value distribution}\label{sec:approx}

This section is devoted to the proof of Theorem \ref{thm:singular-value}.
Our main object to work on is the singular value measure $\nu_{{M}-z}$. Recall that $H(z)$ is the $2n \times 2n$ Hermitian matrix
\begin{equation*}
H(z) := \begin{bmatrix} 0 &
M-z \\  (M-z)^* & 0 \end{bmatrix}.
\end{equation*}
We denote the resolvent $G_n$ as
\begin{equation}\label{eqn:resolvent-original}
G_n(\eta)=(H(z)-\eta)^{-1}, \qquad\eta\in\mathbb{C}^+. 
\end{equation}
The study of $\nu_{{M}-z}$ focuses on the the Stieltjes transform $S_{\widetilde{\nu}_{M-z}}$ of $\widetilde{\nu}_{M-z}$ defined in \eqref{eqn:intro-resolvent-main}, which can be rewritten as
\begin{equation}
  \label{defn:S-transform-singular-V1}
	S_{\widetilde{\nu}_{M-z}}(\eta)=\frac{1}{2n}\text{Tr}(G_n (\eta)), \qquad\eta\in\mathbb{C}^+. 
\end{equation}

It is more convenient to work on block matrices. To this end, we
define the set $\mathbb{H}_+$ as 
\begin{equation*}
\mathbb{H}_+ := \left\{  \begin{bmatrix} \eta &
z \\  \bar{z} & \eta \end{bmatrix} : z \in \mathbb{C}, \eta \in \mathbb{C}_+ \right\}.
\end{equation*}
For $q \in \mathbb{H}_+,$ with
\begin{equation*}
q(z, \eta) := \begin{bmatrix} \eta &
z \\  \bar{z} & \eta \end{bmatrix},
\end{equation*} 
we consider the matrix 
$$B- \un_n \otimes q(z, \eta) \in \mathcal{M}_n (\mathcal{M}_2 (\mathbb{C})),$$
where $B$ is obtained from the $2\times 2$ blocks $B_{ij}$ ($1\leq i, j \leq n$) defined as
\begin{equation*}
B_{ij} := \begin{bmatrix} 0 &
M_{ij} \\  \bar{M}_{ji} & 0 \end{bmatrix}.
\end{equation*}
By permutating the entries, the matrix $H(z) - \eta$ is equivalent to the matrix $B- \un_n \otimes q(z, \eta).$
Hence, instead of studying the original resolvent $G_n$, it is equivalent to consider the following resolvent of $B- \un_n \otimes q(z, \eta)$ defined as
\begin{equation}\label{eq:resolvent}
R(q) = (B- \un_n \otimes q(z, \eta))^{-1}.
\end{equation}
In this paper, the resolvent $R(q)$ is always considered as a $n \times n$ block matrix whose entries are $2 \times 2$ matrices. 
We have the following deterministic bound for the operator norm of the resolvent
\begin{equation} 
\Vert R(q) \Vert \leq \frac{1}{ \Im (\eta)}, \qquad q\in\mathbb{H}^+. 
\end{equation}

We will frequently use the resolvent identity
\begin{equation}
A^{-1} -B^{-1} = A^{-1} (B-A) B^{-1}
\end{equation}
for invertible matrices $A$ and $B$.
We will also use the following Schur’s complement formula for block matrices frequently:
\begin{equation}
\begin{bmatrix} A &
B \\  C & D \end{bmatrix}^{-1} = \begin{bmatrix} A^{-1} + A^{-1} B (D - CA^{-1} B)^{-1} C A^{-1} &
-A^{-1} B (D-CA^{-1}B)^{-1} \\  -(D-CA^{-1}B)^{-1} CA^{-1} & (D - CA^{-1} B)^{-1} \end{bmatrix}
\end{equation}
where $A, B, C, D$ are matrix sub-blocks and $A, D - CA^{-1} B$ are invertible.

\subsection{Approximate fixed point equation}\label{sec:approx-fixed-pt}

For $q \in \mathbb{H}_+,$ let $R(q)_{kk}$ be the $k$-th diagonal $2 \times 2$ block of $R(q),$ and set 
\begin{equation}\label{eq:R_kk}
R(q)_{kk} = \begin{bmatrix} a_k(q) &
b_k(q) \\  c_k(q) & d_k(q) \end{bmatrix} \in \mathcal{M}_2(\mathbb{C}).
\end{equation}
By \eqref{defn:S-transform-singular-V1} and \eqref{eq:resolvent}, we deduce that
\[
  S_{\widetilde{\nu}_{M-z}}(\eta)=\frac{1}{2n}\text{Tr}(R(q))=\frac{1}{2n}\left(\sum_{k=1}^n a_k(q)+\sum_{k=1}^n d_k(q) \right).
\]
One can check that (see \cite[Lemma 4.19]{Bordenave-Chafai-circular} for example)
\begin{equation}\label{eq:entries-resovlent}
\begin{split}
a(q) & := \frac{1}{n} \sum_{k=1}^n a_k(q) = \frac{1}{n} \sum_{k=1}^n d_k(q),\\
b(q) & := \frac{1}{n} \sum_{k=1}^n b_k(q) = \frac{1}{n} \sum_{k=1}^n \bar{c}_k(q),
\end{split}
\end{equation}
where $a(q), b(q)$ denote the common sums. It follows that
\begin{equation}\label{defn:S-transform-singular-V2}
 S_{\widetilde{\nu}_{M-z}}(\eta)=a(q).
\end{equation}
We denote 
\begin{equation}\label{eq:resolvent-1}
\Gamma_n (q) := \frac{1}{n} \sum_{k=1}^n R_{kk}(q)= \begin{bmatrix} a(q) &
b(q) \\  \bar{b}(q) & a(q) \end{bmatrix}.
\end{equation}

The first step towards the proof of Theorem \ref{thm:singular-value} is to obtain the following fixed point equation for the Stieltjes transform of $H(z)$. 

\begin{theorem}\label{thm:fix-point-equ} 
For any $q=q(z,\eta)\in \mathbb{H}_+,$ let $a(q), b(q)$ be given in \eqref{eq:entries-resovlent}, then any accumulation point of 
\begin{equation*}
\mathbb{E} \begin{bmatrix} a(q) &
b(q) \\  \bar{b}(q) & a(q) \end{bmatrix} 
\end{equation*}
is a solution of the fixed point equation
\begin{equation}\label{eq:fixed-point-eq}
\begin{bmatrix} \alpha &
\beta \\  \bar{\beta} & \alpha \end{bmatrix} = - \mathbb{E} \left( q+ \begin{bmatrix} \alpha &
\gamma \bar{\beta} \\   \gamma \beta & \alpha \end{bmatrix}- \begin{bmatrix} 0 &
G \\ \overline{G} & 0 \end{bmatrix} \right)^{-1},
\end{equation}
where $G$ is a Gaussian random variable on $\mathbb{C}$ with law $\mathcal{N}(0, K_\mu)$ and covariance matrix $K_\mu=\begin{bmatrix} \mu &
0 \\  0 & 1-\mu \end{bmatrix}.$ 
\end{theorem}

The general approach for the proof of Theorem \ref{thm:fix-point-equ} is similar to the proof for \cite[Section 2.5]{BordenaveCaputoChafai2014markov}. 
However, we need extra work to deal with the dependency of the entries, and thus the proof is substantially more technical. For $\vert \gamma \vert \leq 1$, we denote $\Sigma$ as an operator on $2\times 2$ matrices given by
	\begin{equation}
	\Sigma \begin{bmatrix} a &
	b \\  c & d \end{bmatrix}: = \begin{bmatrix} d &
	\gamma c \\  \gamma b & a \end{bmatrix}.
	\end{equation} 
We will frequently use the fact that $\Vert \Sigma\Vert\leq C$ for some constant $C$, where $\Vert \Sigma\Vert$ denotes the operator norm of $\Sigma$.  
Our goal is to show that
\begin{equation}\label{eqn:goal-limit1}
 \mathbb{E} [\Gamma_n(q)]  = - \mathbb{E}  \left( q + \mathbb{E} \left[ \Sigma(\Gamma_n(q)) \right] -  \begin{bmatrix} 0 &
G \\  \overline{G} & 0 \end{bmatrix}  \right)^{-1} + \varepsilon
\end{equation}
with $\Vert \varepsilon \Vert = o(1).$ Recall that $\Gamma_n (q) = \frac{1}{n} \sum_{k=1}^n R(q)_{kk},$ it suffices to obtain a similar estimation for $R(q)_{kk}$ for any $q=q(z,\eta)\in \mathbb{H}_+.$ 

For $k=1, \ldots, n$, denote $Q^{(k)} \in \mathcal{M}_{n-1,1}(\mathcal{M}_2(\mathbb{C}))$ by the $2(n-1) \times 2$ matrix whose blocks are given by
\begin{equation}\label{eqn:defn-Qk}
Q_i^{(k)} = \begin{bmatrix} 0 &
M_{ki} \\  \bar{M}_{ik} & 0 \end{bmatrix}  =  \frac{1}{\sqrt{n}} \begin{bmatrix} 0 &
X_{ki} \\  \bar{X}_{ik} & 0 \end{bmatrix}, i=1, 2, \ldots, k-1, k+1, \ldots, n.
\end{equation}
Let $B^{(k)}$ be the minor of $B$, where the $k^{th}$ row and $k^{th}$ column of $B$ have been removed, and $R^{(k)} (q) = (B^{(k)} - \un_{n-1} \otimes q)^{-1}$ be the resolvent of the minor. 
By the Schur block inversion formula, for any $k=1, \ldots, n,$ we have 
\begin{equation}\nonumber
R(q)_{kk} = \left( \begin{bmatrix} 0 &
M_{kk} \\  \bar{M}_{kk} & 0 \end{bmatrix} -q - {Q^{(k)}}^* R^{(k)}(q) Q^{(k)} \right)^{-1}.
\end{equation}

However, the entries of $Q^{(k)}$ and $R^{(k)}$ have some correlations which prevents us from working on them directly. We instead introduce below some modifications to help us prove \eqref{eqn:goal-limit1}. 
We denote 
\begin{equation}\label{eq:resolvent-2}
\Gamma_n^{(k)} (q) := \frac{1}{n} \sum_{i\neq k} R_{ii}^{(k)} (q).
\end{equation}
Define the matrix $\widetilde{M}^{(k)} \in \mathcal{M}_{n-1}(\mathbb{C})$ by 
$$\widetilde{M}^{(k)}_{ij} = \frac{X_{ij}}{\sqrt{n}} -\delta_{i,j} \sum_{l \neq k} \frac{X_{il}}{ \sqrt{n}}$$
for $i,j=1, \ldots, k-1, k+1, \ldots, n.$ Let $\widetilde{B}^{(k)}$ and $\widetilde{R}^{(k)}$ in $\mathcal{M}_{n-1}(\mathcal{M}_2(\mathbb{C}))$ be the matrices similarly obtained by replacing $M$ with $\widetilde{M}^{(k)}.$ Denote
\begin{equation}\label{eq:resolvent-3}
\widetilde{\Gamma}^{(k)}_n(q) : = \frac{1}{n} \sum_{i\neq k} \widetilde{R}^{(k)}_{ii} (q).
\end{equation}

Our proof is mainly divided into two parts. In the first part, by obtaining some a priori bounds for the resolvents, we can show that 
$$\left\Vert \widetilde{\Gamma}_n^{(k)} (q)-  \Gamma_n(q) \right\Vert = o(1), \; \text{in probability}$$
for any $k=1, \ldots, n$ (see Lemma \ref{lem:epsilon2} and \ref{lem:epsilon3}). Then Lemma \ref{lem:epsilon5} says that
$\Gamma_n (q)$ concentrates around its expectation $\mathbb{E} [\Gamma_n(q)].$ The second part devotes to estimating the difference between ${Q^{(k)}}^* \widetilde{R}^{(k)} Q^{(k)},$
and ${Q^{(k)}}^* R^{(k)} Q^{(k)}$ with respect to the operator norm. By a careful analysis on the entries (see Lemma \ref{lem:epsilon4}), we have 
$$\left\Vert {Q^{(k)}}^* \widetilde{R}^{(k)} Q^{(k)} - \Sigma (\widetilde{\Gamma}_n^{(k)}(q)) \right\Vert = o(1), \; \text{in probability}$$
for any $k=1, \ldots, n.$
Hence by combining the first part, we can show that 
$$\left\Vert {Q^{(k)}}^* R^{(k)} Q^{(k)} - \mathbb{E} \left[ \Sigma(\Gamma_n(q)) \right] \right\Vert = o(1), \; \text{in probability}$$
for any $k=1, \ldots, n.$ In addition, by the central limit theorem, $M_{kk}$ converges to $G$ as $n \rightarrow \infty$. Equation \eqref{eqn:goal-limit1} would follow from these results. 


\subsubsection{Bounds for the resolvents of the Laplacian and its minor}\label{subsubsec:minor}

In this subsection, we estimate differences between $\Gamma_n(q),$ $\Gamma_n^{(k)} (q)$, and $\widetilde{\Gamma}^{(k)}_n(q)$ given by Equation \eqref{eq:resolvent-1}, \eqref{eq:resolvent-2}, and \eqref{eq:resolvent-3} respectively. 

\begin{lemma}\label{lem:epsilon2}
With probability $1$, we have a priori bound 
\begin{equation}\label{eq:norm-epsilon2}
\sup_{1 \leq k \leq n} \left\Vert \widetilde{\Gamma}_n^{(k)} (q) - \Gamma^{(k)}_n(q) \right\Vert = O\left( \frac{\kappa(n)}{\sqrt{n}\Im (\eta)^{2}} \right).
\end{equation}
\end{lemma}

\begin{proof}
By the resolvent identity we have 
\begin{equation*}
\begin{split}
\Vert  \widetilde{\Gamma}_n^{(k)} (q) - \Gamma^{(k)}_n(q)  \Vert  & = \frac{1}{n} \left\Vert \sum_{i\neq k} \widetilde{R}^{(k)}_{ii} - \sum_{i \neq k} R^{(k)}_{ii}  \right\Vert \leq \Vert \widetilde{R}^{(k)}-R^{(k)}\Vert \\
& = \Vert \widetilde{R}^{(k)}( B^{(k)} - \widetilde{B}^{(k)}) R^{(k)}\Vert \\
& \leq  \frac{1}{\sqrt{n} \Im (\eta)^2} \cdot \left\Vert {\rm diag} \; \left(  \begin{bmatrix} 0 &
X_{ki} \\  \bar{X}_{ik} & 0 \end{bmatrix} \right)_{i \neq k}\right\Vert \\
&  = \frac{\kappa(n)}{\sqrt{n} \Im (\eta)^2}. 
\end{split}
\end{equation*}
\end{proof}

The following inequality is useful when we compare $\Gamma^{(k)}_n(q)$ with  $\Gamma_n(q)$. 
\begin{lemma}\label{lem:pertubation-resolvent}
Let $A, B$ be two $2n \times 2n$ Hermitian matrices. If $q = \begin{bmatrix} \eta &
z \\  \bar{z} & \eta \end{bmatrix} \in \mathbb{H}_{+},$ $R_A(q)= (A- \un_n \otimes q)^{-1}$ and $R_B(q) = (B- \un_n \otimes q)^{-1},$ then we have
\begin{equation}
\sum_{k=1}^n \left\Vert R_A(q)_{kk} - R_B(q)_{kk} \right\Vert \leq 2 \frac{{\rm rank} \; (A-B)}{\Im (\eta)}.
\end{equation}
\end{lemma}

\begin{proof}
By the resolvent identity, we have
\begin{equation*}
M:= R_A - R_B = R_A(B-A) R_B.
\end{equation*}
It follows that $r= {\rm rank} \; (M) \leq {\rm rank} \; (B-A).$ We also note that $\Vert M \Vert \leq 2 \Im (\eta)^{-1}.$ Hence, by the singular value decomposition of $M,$ we have 
$$M = \sum_{i=1}^r s_i u_i v_i^*,$$
where $s_1, \ldots, s_r$ are singular values of $M,$ and $u_1, \ldots, u_r$ and $v_1, \ldots, v_r$ are the associated orthonormal left or right singular vectors. 

Hence 
$$R_A(q)_{kk} - R_B(q)_{kk} = M_{kk} =  \sum_{i=1}^r s_i  e_k^*u_i v_i^*e_k,$$
where $e_k$ is the $n \times 1$ vector whose $k^{th}$ $2\times 2$ block is the identity matrix and whose other entries are zero. 
We obtain from the Cauchy-Schwarz inequality, 
\begin{equation*}
\begin{split}
\sum_{k=1}^n \Vert M_{kk} \Vert & \leq \Vert M \Vert \sum_{i=1}^r \sum_{k=1}^n \Vert e_k^*u_i v_i^*e_k \Vert \\
& \leq \Vert M \Vert \sum_{i=1}^r \sqrt{\sum_{k=1}^n \Vert e_k^*u_i\Vert^2} \cdot \sqrt{\sum_{k=1}^n \Vert v_i^*e_k \Vert^2} \\
& = \Vert M \Vert \sum_{i=1}^r \sqrt{\sum_{k=1}^n u_i^* e_k e_k^*u_i} \cdot \sqrt{\sum_{k=1}^n v_i^*e_k e_k^* v_i} \\
& = r \Vert M \Vert. 
\end{split}
\end{equation*}
\end{proof}

\begin{lemma}\label{lem:epsilon3}
With probability $1$, we have a priori bound 
\begin{equation}\label{eq:norm-epsilon3}
\sup_{1 \leq k \leq n} \Vert \Gamma^{(k)}_n(q) - \Gamma_n(q) \Vert = O(n^{-1}\Im (\eta)^{-1}).
\end{equation}
\end{lemma}

\begin{proof}
We follow the idea of \cite[Lemma 6.4]{2013Low}. We write $B_k$ to be the $k^{th}$ column (of $2 \times 2$ blocks) of $B$ and define the modified resolvent 
\begin{equation}
\check{R}^{(k)}(q) : =   ( (B -e_k B_k^*-B_k e_k^*) - \un_n \otimes q)^{-1},
\end{equation}
where $e_k$ is the $k \times 1$ vector whose $k^{th}$ $2\times 2$ block is the identity matrix and whose other entries are zero. 
Hence,
$$\left\Vert \sum_{i\neq k} R^{(k)}_{ii} - \sum_{i=1}^{n} \check{R}_{ii}^{(k)} \right\Vert  = \Vert (B_{kk}+q)^{-1} \Vert \leq \frac{1}{\Im (\eta)}.$$

Moreover, by the resolvent identity, 
\begin{equation*}
\begin{split}
R(q) - \check{R}^{(k)}(q)  =  R(q) (e_k B_k^* + B_k e_k^*) \check{R}^{(k)}(q),
\end{split}
\end{equation*}
The matrix $R(q) - \check{R}^{(k)} (q)$ has rank at most $2$ (viewed as a $2n \times 2n$ matrix). Using Lemma \ref{lem:pertubation-resolvent},
$$\sum_{i=1}^n \left\Vert \check{R}^{(k)}_{ii} - R_{ii} \right\Vert \leq \frac{4} {\Im (\eta)}.$$
Thus, we obtain the estimate
\begin{equation*}
\begin{split}
\Vert \Gamma^{(k)}_n(q) - \Gamma_n(q) \Vert & = \frac{1}{n} \left\Vert \sum_{i\neq k} R^{(k)}_{ii} - \sum_{i=1}^n R_{ii} \right\Vert \\
& \leq \frac{1}{n} \left\Vert \left( \sum_{i\neq k} R^{(k)}_{ii} - \sum_{i=1}^n \check{R}_{ii}^{(k)} \right) + \left(  \sum_{i=1}^n \check{R}_{ii}^{(k)} - \sum_{i=1}^{n} R_{ii} \right) \right\Vert \\
& \leq \frac{1}{n} \frac{5} {\Im (\eta)}.
\end{split}
\end{equation*}
\end{proof}

\begin{lemma}\label{lem:epsilon5}
For any $q=q(z,\eta)\in \mathbb{H}_+,$ the norm $\left\Vert \Gamma_n (q) - \mathbb{E} \left[ \Gamma_n (q)\right] \right\Vert$ converges to 0 in probability as $n \rightarrow \infty$.
\end{lemma}

\begin{proof}
We adapt the methods in \cite[Lemma 7.7]{NguyenORourke2015} and \cite[Lemma A. 2]{Tikhom2023}. Denote $\mathbb{E}_k$ by the conditional expectation with respect to the $\sigma$-algebra generated by $\{X_{ij}: 1\leq i, j \leq k\}.$ Now we decompose $\Gamma_n (q) - \mathbb{E} \left[ \Gamma_n (q)\right]$ into the sum of martingale differences:
\begin{equation*}
\begin{split}
\Gamma_n (q) - \mathbb{E} \left[ \Gamma_n (q) \right] & = \sum_{k=1}^n (\mathbb{E}_k - \mathbb{E}_{k-1} )  [\Gamma_n(q)]\\
& = \sum_{k=1}^n (\mathbb{E}_k - \mathbb{E}_{k-1} )  [\Gamma_n(q) - \widetilde{\Gamma}^{(k)}_n(q)]\\
& = \sum_{k=1}^n (\mathbb{E}_k - \mathbb{E}_{k-1} )  [\Gamma_n - \Gamma^{(k)}_n + \Gamma^{(k)}_n- \widetilde{\Gamma}^{(k)}_n],
\end{split}
\end{equation*}
where we used the equality 
\begin{equation*}
\mathbb{E}_k \left[ \widetilde{\Gamma}_n^{(k)}(q) \right] = \mathbb{E}_{k-1} \left[ \widetilde{\Gamma}^{(k)}_n(q) \right].
\end{equation*} 
Hence, denote $\gamma_k^{(1)} : = (\mathbb{E}_k - \mathbb{E}_{k-1} )  [\Gamma_n - \Gamma^{(k)}_n]$ and $\gamma_k^{(2)} : = (\mathbb{E}_k - \mathbb{E}_{k-1} )  [\Gamma^{(k)}_n- \widetilde{\Gamma}^{(k)}_n],$ it suffices to show that $\Vert \sum_{k=1}^n \gamma^{(i)}_k \Vert$ converges to $0$ in probability, $i=1,2$.

For any $2 \times 2$ matrix $A$, denote $A^{ab}$ by the ${ab}^{th}$ entry of $A$, $a, b =1, 2.$
By the Burkholder inequality, there exists an universal constant $C_0>0$ such that 
\begin{equation*}
\begin{split}
\mathbb{E} & \left[  \vert \sum_{k=1}^n \gamma^{(1) ab}_k  \vert^4 \right]  = \mathbb{E} \left[  \left\vert \sum_{k=1}^n (\mathbb{E}_{k} - \mathbb{E}_{k-1}) \left[ \Gamma_n^{ab} -\Gamma_n^{(k)ab} \right] \right\vert^4 \right] \\
& \leq C_0 \cdot \mathbb{E} \left[  \sum_{k=1}^n  \left\vert (\mathbb{E}_k- \mathbb{E}_{k-1}) \left[ \Gamma_n^{ab} -\Gamma_n^{(k) ab} \right] \right\vert^2 \right]^2.
\end{split}
\end{equation*}
Note that $\Vert \Gamma^{(k)}_n(q) - \Gamma_n(q) \Vert = O(n^{-1}\Im (\eta)^{-1})$ (see Lemma \ref{lem:epsilon3}), it follows that
\begin{equation*}
\begin{split}
\left\vert (\mathbb{E}_k - \mathbb{E}_{k-1}) \left[ \Gamma_n^{ab} -\Gamma_n^{(k) ab} \right] \right\vert \leq \Vert \Gamma^{(k)}_n(q) - \Gamma_n(q) \Vert = O(n^{-1}\Im (\eta)^{-1}).
\end{split}
\end{equation*}
We hence conclude that for any $a, b=1,2$,
\begin{equation*}
\mathbb{E} \left[  \vert  \sum_{k=1}^n \gamma^{(1) ab}_k \vert^4 \right] \leq C_0 \cdot O(n^{-2} \Im (\eta)^{-4}).
\end{equation*}
By the Markov inequality and Borel-Cantelli lemma, almost surely as $n \rightarrow \infty,$
\begin{equation}\label{eq:estimate-1}
\left\Vert  \sum_{k=1}^n \gamma^{(1)}_k  \right\Vert \rightarrow 0.
\end{equation}

Now we move to bound $\Vert \sum_{k=1}^n \gamma^{(2)}_k \Vert.$ It suffices to bound 
$\sum_{k=1}^n \mathbb{E} [\vert \gamma^{(2)ab}_k \vert^2].$ 
Note that 
\begin{equation*}
\mathbb{E} [\vert \gamma^{(2)ab}_k \vert^2] \leq 2 \mathbb{E} [ \vert\Gamma_n^{(k) ab} -\widetilde{\Gamma}_n^{(k)ab} \vert^2 ].
\end{equation*}
Denote an $2(n-1) \times 2(n-1)$ diagonal matrix $D^{(k)}$ with diagonal entries
\begin{equation*}
D^{(k)} = B^{(k)} - \widetilde{B}^{(k)} = {\rm diag} \; \left(  \begin{bmatrix} 0 &
X_{ki}/\sqrt{n} \\  \bar{X}_{ik}/ \sqrt{n} & 0 \end{bmatrix} \right)_{i \neq k}.
\end{equation*}
We have 
\begin{equation*}
\begin{split}
\widetilde{\Gamma}_n^{(k)} - \Gamma^{(k)}_n & = \frac{1}{n} \sum_{i \neq k} (\widetilde{R}^{(k)}_{ii} - R^{(k)}_{ii} ) = \frac{1}{n} \sum_{i \neq k} (\widetilde{R}^{(k)} D^{(k)} R^{(k)})_{ii}\\
& = \frac{1}{n} \sum_{i \neq k} (\widetilde{R}^{(k)} D^{(k)} \widetilde{R}^{(k)})_{ii} - \frac{1}{n} \sum_{i \neq k} (\widetilde{R}^{(k)} D^{(k)} R^{(k)} D^{(k)} \widetilde{R}^{(k)})_{ii} \\
& : = \frac{1}{n} \sum_{i \neq k} A^{(k)}_{ii}  - \frac{1}{n} \sum_{i \neq k} B^{(k)}_{ii}.
\end{split}
\end{equation*}
This implies that 
\begin{equation}\label{eq:bound-1}
\mathbb{E} [\vert \gamma^{(2)ab}_k \vert^2] \leq \frac{4}{n^2} \mathbb{E} \left[ \vert \sum_{i\neq k} A^{(k) ab}_{ii} \vert^2 +\vert \sum_{i\neq k} B^{(k)ab}_{ii} \vert^2  \right].
\end{equation}

Note that 
\begin{equation*}
\begin{split}
A^{(k) ab}_{ii} & =  \sum_{j \neq k}  \sum_{c, d=1,2} \widetilde{R}^{(k)ac}_{ij} D^{(k)cd}_{jj} \widetilde{R}^{(k)da}_{ji}\\
& = \frac{1}{\sqrt{n}} \sum_{j \neq k} (X_{kj} \widetilde{R}^{(k)a1}_{ij} \widetilde{R}^{(k)2a}_{ji} + \bar{X}_{jk} \widetilde{R}^{(k)a2}_{ij} \widetilde{R}^{(k)1a}_{ji}).
\end{split}
\end{equation*}
It is straightforward to check that 
\begin{equation}
\vert \sum_{i \neq k}  \widetilde{R}^{(k)a1}_{ij} \widetilde{R}^{(k)2a}_{ji} \vert^2  \leq \Vert (\widetilde{R}^{(k)} \widetilde{R}^{(k)})^2 \Vert \leq \frac{1}{\Im (\eta)^4},
\end{equation}
and similarly $\vert \sum_{i \neq k} \widetilde{R}^{(k)a2}_{ij} \widetilde{R}^{(k)1a}_{ji} \vert^2 \leq 1/ \Im (\eta)^4.$

Therefore, by the independence of random variables $X_{kj}$ (resp. $X_{jk}$) for $j \neq k$ and $\widetilde{R}^{(k)}$, we have 
\begin{equation*}
\begin{split}
\mathbb{E}  \vert \sum_{i \neq k} A_{ii}^{(k) ab} \vert^2 & \leq \frac{2}{n} \mathbb{E} \left[ \vert \sum_{i, j \neq k} X_{kj} \widetilde{R}^{(k)a1}_{ij} \widetilde{R}^{(k)2a}_{ji} \vert^2 + \vert \sum_{i, j \neq k} \bar{X}_{jk} \widetilde{R}^{(k)a2}_{ij} \widetilde{R}^{(k)1a}_{ji}\vert^2 \right] \\
& \leq \frac{2}{n} \cdot \frac{1}{\Im (\eta)^4} \sum_{j \neq k} \mathbb{E} [\vert X_{kj} \vert^2 + \vert \bar{X}_{jk} \vert^2] \leq \frac{4}{\Im (\eta)^4}.
\end{split}
\end{equation*}

For the second term of \eqref{eq:bound-1}, we have
\begin{equation*}
\begin{split}
B^{(k) ab }_{ii} & =  \sum_{j, p\neq k}  \sum_{c, d,e, f=1,2} \widetilde{R}^{(k)ac}_{ij} D^{(k)cd}_{jj} R^{(k)de}_{jp} D^{(k)ef}_{pp} \widetilde{R}^{(k)fb}_{pi} \\
& = \frac{1}{n} \sum_{j, p\neq k}  (  X_{kj} X_{kp} \widetilde{R}^{(k)a1}_{ij} R^{(k)21}_{jp} \widetilde{R}^{(k)2b}_{pi} + 
X_{kj} \bar{X}_{pk} \widetilde{R}^{(k)a1}_{ij} R^{(k)22}_{jp} \widetilde{R}^{(k)1b}_{pi}) \\
& +  \frac{1}{n} \sum_{j, p\neq k} ( \bar{X}_{jk} \bar{X}_{pk}  \widetilde{R}^{(k)a2}_{ij} R^{(k)12}_{jp} \widetilde{R}^{(k)1b}_{pi}+
 \bar{X}_{jk}  X_{kp} \widetilde{R}^{(k)a2}_{ij} R^{(k)11}_{jp}  \widetilde{R}^{(k)2b}_{pi}).
\end{split}
\end{equation*}
Firstly, by using the Cauchy-Schwarz inequality, we have
\begin{equation*}
\begin{split}
\frac{1}{n^2} \mathbb{E} & \left[ \vert \sum_{i, j, p\neq k}  X_{kj} X_{kp} \widetilde{R}^{(k)a1}_{ij} R^{(k)21}_{jp} \widetilde{R}^{(k)2b}_{pi} \vert^2 \right] \\
& \leq \frac{1}{n} \mathbb{E} \left[ \sum_{j\neq k} ( \vert X_{kj} \vert^2 \cdot \vert \sum_{i, p \neq k} X_{kp} \widetilde{R}^{(k)a1}_{ij} R^{(k)21}_{jp} \widetilde{R}^{(k)2b}_{pi} \vert^2) \right] \\
& \leq \frac{\kappa^2(n)}{n} \mathbb{E} \left[ \sum_{j\neq k}  (\vert X_{kj} \vert^2 \cdot \sum_{p \neq k} \vert R^{(k)21}_{jp}  \vert^2 \cdot \sum_{p \neq k} (\sum_{i \neq k}\widetilde{R}^{(k)a1}_{ij} \widetilde{R}^{(k)2b}_{pi}  )^2 ) \right].
\end{split}
\end{equation*}

Denote $R^{(k)ab}$ by the $n \times n$ matrix formed from taking each $R^{(k)}_{ij}$ block and replace it by its $ab$-th entry (see \cite[Section 6.1]{2013Low}), then we have
\begin{equation}
\begin{split}
\sum_{p \neq k} \vert R^{(k)21}_{jp}  \vert^2 & = (R^{(k)21} {R^{(k)21}}^*)_{jj} \leq \Vert R^{(k)21} {R^{(k)21}}^* \Vert \\
& \leq \Vert R^{(k)} {R^{(k)}}^* \Vert \leq \frac{1}{\Im (\eta)^2},
\end{split}
\end{equation}
and similarly 
\begin{equation}
\sum_{p \neq k} (\sum_{i \neq k}\widetilde{R}^{(k)a1}_{ij} \widetilde{R}^{(k)2b}_{pi}  )^2   \leq  \Vert (\widetilde{R}^{(k)} {\widetilde{R}^{(k)}}) \cdot (\widetilde{R}^{(k)} {\widetilde{R}^{(k)}})^* \Vert \leq \frac{1}{\Im (\eta)^4}.
\end{equation}
It follows that 
\begin{equation}
\frac{1}{n^2} \mathbb{E}  \left[ \vert \sum_{i, j, p\neq k}  X_{kj} X_{kp} \widetilde{R}^{(k)a1}_{ij} R^{(k)21}_{jp} \widetilde{R}^{(k)2b}_{pi} \vert^2 \right] \leq \frac{\kappa^2(n) }{\Im (\eta)^6}. 
\end{equation}
We can similarly bound the rest of the terms of $\sum_{i\neq k} B^{(k) ab }_{ii}.$ Thus, we obtain 
\begin{equation}
\mathbb{E}  \vert \sum_{i \neq k} B_{ii}^{(k) ab} \vert^2 \leq \frac{16 \kappa^2(n) }{\Im (\eta)^6}.
\end{equation}

In summary, we conclude that for any $a, b =1,2$,
\begin{equation}
\begin{split}
\sum_{k=1}^n \mathbb{E} [\vert \gamma^{(2)ab}_k \vert^2] & \leq \frac{4}{n^2} \sum_{k=1}^n \mathbb{E} \left[ \vert \sum_{i\neq k} A^{(k) ab}_{ii} \vert^2 +\vert \sum_{i\neq k} B^{(k)ab}_{ii} \vert^2  \right]\\
& \leq \frac{4}{n} \left(\frac{4}{\Im (\eta)^4} + \frac{16 \kappa^2(n) }{\Im (\eta)^6} \right) \rightarrow 0, \; \text{as} \; n \rightarrow \infty.
\end{split}
\end{equation}
By the Markov inequality, the norm $\Vert  \sum_{k=1}^n \gamma^{(2)}_k \Vert$ converges to $0$ in probability as $n \rightarrow \infty.$ Combining \eqref{eq:estimate-1}, we complete our proof.
\end{proof}

\subsubsection{Approximate fixed point equation}

In this subsection, we give the proof for Theorem \ref{thm:fix-point-equ}.

\begin{lemma}\label{lem:epsilon4}
Let $1 \leq k \leq n$, and fix $q=q(z,\eta)\in\mathbb{H}_+$.  Recall that $Q^{(k)}$ is defined in \eqref{eqn:defn-Qk}. 
We denote by $\mathcal{F}_{k}$ the $\sigma$-algebra spanned by the variables $\{ X_{ij}: i, j = 1, 2, \ldots, k-1, k+1, \ldots, n\}, $ and $\mathbb{E}_{k} [\; \cdot \; ] : = \mathbb{E}[\; \cdot \; |\mathcal{F}_{k}].$ We have 
\begin{enumerate}[\rm (1)]
\item  $\mathbb{E}_{k} \left[ {Q^{(k)}}^* \widetilde{R}^{(k)} Q^{(k)} \right]  = \Sigma (\widetilde{\Gamma}_n^{(k)}(q)).$
\item  $\left\Vert {Q^{(k)}}^* \widetilde{R}^{(k)} Q^{(k)} - \mathbb{E}_{k} \left[ {Q^{(k)}}^* \widetilde{R}^{(k)} Q^{(k)} \right] \right\Vert\rightarrow 0$ in probability as $n \rightarrow \infty.$
\end{enumerate}
\end{lemma}

\begin{proof}
{\rm (1)} Let $$ \widetilde{R}_{ij}^{(k)} =  \begin{bmatrix} \tilde{a}_{ij} &
\tilde{b}_{ij} \\  \tilde{c}_{ij} & \tilde{d}_{ij} \end{bmatrix}, \; i, j =1, \ldots, k-1, k+1, \ldots, n,$$
then we have
\begin{equation*}
\begin{split}
\mathbb{E}_{k} \left[ {Q^{(k)}}^* \widetilde{R}^{(k)} Q^{(k)} \right] & = \sum_{i,j\neq k} \mathbb{E}_{k} \left[ {Q^{(k)}_i}^* \widetilde{R}^{(k)}_{ij} Q^{(k)}_j \right] \\
& = \frac{1}{n} \sum_{i,j \neq k} \mathbb{E}_{k} \begin{bmatrix} \tilde{d}_{ij} X_{ik} \bar{X}_{jk} &
\tilde{c}_{ij} X_{ik} X_{kj} \\  \tilde{b}_{ij} \bar{X}_{ki} \bar{X}_{jk} & \tilde{a}_{ij} \bar{X}_{ki} X_{kj} \end{bmatrix}\\
& =  \frac{1}{n} \sum_{i\neq k} \begin{bmatrix} \tilde{d}_{ii}  &
\gamma \tilde{c}_{ii}  \\  \gamma \tilde{b}_{ii}  & \tilde{a}_{ii} \end{bmatrix}\\ 
& =  \Sigma (\widetilde{\Gamma}_n^{(k)}(q)).
\end{split}
\end{equation*}

\noindent {\rm (2)} Denote 
$$Z_k : = \mathbb{E}_{k} \left\vert {Q^{(k)}}^* \widetilde{R}^{(k)} Q^{(k)} - \mathbb{E}_{k} \left[ {Q^{(k)}}^* \widetilde{R}^{(k)} Q^{(k)} \right] \right\vert^2,$$ 
it suffices to show that $\Vert Z_k \Vert$ converges to $0$ with probability $1,$ as $n \rightarrow \infty.$ By an elementary computation, we have 
\begin{equation*}
\begin{split}
 {\rm Tr} (Z_k) & = \frac{1}{n^2} \mathbb{E}_{k}  \Big[ \vert \sum_{i,j \neq k} \tilde{d}_{ij} X_{ik} \bar{X}_{jk}  - \sum_{i \neq k} \tilde{d}_{ii} \vert^2  + \vert \sum_{i,j \neq k}  \tilde{b}_{ij} \bar{X}_{ki} X_{jk}  - \gamma \sum_{i \neq k} \tilde{b}_{ii} \vert^2  \\
 & + \vert \sum_{i,j \neq k}  \tilde{c}_{ij} \bar{X}_{ik} \bar{X}_{kj} - \gamma \sum_{i \neq k} \tilde{c}_{ii} \vert^2  +  \vert \sum_{i,j \neq k} \tilde{a}_{ij} \bar{X}_{ki} X_{kj}  - \sum_{i \neq k} \tilde{a}_{ii} \vert^2 \Big]\\
 & : = (\mathrm{I}) + (\mathrm{II}) + (\mathrm{III}) + (\mathrm{IV}).
\end{split}
\end{equation*}

First, we can expand term ($\mathrm{I}$) as
\begin{equation*}
\begin{split}
(\mathrm{I})  = \frac{1}{n^2} \left( \sum_{i,j,p,q \neq k} \mathbb{E} [\bar{X}_{ik} X_{jk} X_{pk} \bar{X}_{qk}] \bar{\tilde{d}}_{ij}  \tilde{d}_{pq}  - \vert \sum_{i\neq k} \tilde{d}_{ii} \vert^2 \right),
\end{split}
\end{equation*}
where $ \bar{\tilde{d}}_{ij}$ is the complex conjugate of $\tilde{d}_{ij}.$ By Condition {\bf C0}, the expectation $\mathbb{E} [\bar{X}_{ik} X_{jk} X_{pk} \bar{X}_{qk}]$ does not vanish whenever 
\begin{enumerate}[\rm (a)]
\item $i=j=p=q,$ $\mathbb{E} [\bar{X}_{ik} X_{jk} X_{pk} \bar{X}_{qk}]= \mathbb{E} [\vert X_{ik} \vert^4] \leq \kappa^2(n);$
\item $i=j \neq p=q,$ $\mathbb{E} [\bar{X}_{ik} X_{jk} X_{pk} \bar{X}_{qk}]= \mathbb{E} [\vert X_{ik} \vert^2] \cdot \mathbb{E} [\vert X_{pk} \vert^2] \leq 1;$
\item $i=p \neq j=q,$ $\mathbb{E} [\bar{X}_{ik} X_{jk} X_{pk} \bar{X}_{qk}]= \mathbb{E} [\vert X_{ik} \vert^2] \cdot \mathbb{E} [ \vert X_{jk} \vert^2] \leq 1;$
\item $i=q \neq j=p,$ $\mathbb{E} [\bar{X}_{ik} X_{jk} X_{pk} \bar{X}_{qk}]= \mathbb{E} [ \bar{X}_{ik}^2] \cdot \mathbb{E} [X_{jk}^2] \leq 1.$
\end{enumerate}
It follows that 
\begin{equation*}
\begin{split}
(\mathrm{I}) & \leq \frac{1}{n^2} \left( \sum_{i \neq j} (\bar{\tilde{d}}_{ii} \tilde{d}_{jj} + \vert \tilde{d}_{ij} \vert^2 + \bar{\tilde{d}}_{ij} \tilde{d}_{ji} )+ \kappa^2(n) \sum_{i\neq k} \vert \tilde{d}_{ii}\vert^2 - \vert \sum_{i\neq k} \tilde{d}_{ii} \vert^2 \right)\\
& = \frac{1}{n^2} \left( \sum_{i \neq j} ( \vert \tilde{d}_{ij} \vert^2 + \bar{\tilde{d}}_{ij} \tilde{d}_{ji} ) + (\kappa^2(n)-1) \sum_{i\neq k} \vert \tilde{d}_{ii}\vert^2 \right).
\end{split}
\end{equation*}

For the term $(\mathrm{II}) $, we have
\begin{equation*}
\begin{split}
(\mathrm{II})  = \frac{1}{n^2} \left( \sum_{i,j,p,q \neq k} \mathbb{E} [X_{ki} \bar{X}_{jk} X_{kp} \bar{X}_{qk}] \bar{\tilde{b}}_{ij} \tilde{b}_{pq}  - \gamma^2 \vert \sum_{i \neq k} \tilde{b}_{ii} \vert^2 \right).
\end{split}
\end{equation*}
Note that the expectation $\mathbb{E} [X_{ki} \bar{X}_{jk} X_{kp} \bar{X}_{qk}]$ does not vanish whenever
\begin{enumerate}[\rm (a)]
\item $i=j=p=q,$ $\mathbb{E} [X_{ki} \bar{X}_{jk} X_{kp} \bar{X}_{qk}]= \mathbb{E} [X_{ki}^2 \bar{X}_{ik}^2] \leq \kappa^2(n);$
\item $i=j \neq p=q,$ $\mathbb{E} [X_{ki} \bar{X}_{jk} X_{kp} \bar{X}_{qk}]= \mathbb{E} [X_{ki}\bar{X}_{ik}] \cdot \mathbb{E} [X_{kp}\bar{X}_{pk}] = \gamma^2;$
\item $i=p \neq j=q,$ $\mathbb{E} [X_{ki} \bar{X}_{jk} X_{kp} \bar{X}_{qk}]= \mathbb{E} [X_{ki}^2] \cdot \mathbb{E} [\bar{X}_{jk}^2] \leq 1;$
\item $i=q \neq j=p,$ $\mathbb{E} [X_{ki} \bar{X}_{jk} X_{kp} \bar{X}_{qk}]= \mathbb{E} [X_{ki}\bar{X}_{ik}] \cdot \mathbb{E} [X_{kj}\bar{X}_{jk}] = \gamma^2.$
\end{enumerate}
It follows that 
\begin{equation*}
\begin{split}
(\mathrm{II}) & \leq \frac{1}{n^2} \left( \sum_{i \neq j} ( \vert \tilde{b}_{ij}\vert^2 + \gamma^2  \bar{\tilde{b}}_{ii} \tilde{b}_{jj} + \gamma^2 \bar{\tilde{b}}_{ij} \tilde{b}_{ji} )+    \sum_{i\neq k} \vert \tilde{b}_{ii} \vert^2 -  \gamma^2 \vert \sum_{i\neq k} \tilde{b}_{ii} \vert^2 \right)\\
& = \frac{1}{n^2} \left( \sum_{i \neq j} ( \vert \tilde{b}_{ij}\vert^2 +\gamma^2 \bar{\tilde{b}}_{ij} \tilde{b}_{ji}) + (\kappa^2(n)- \gamma^2) \sum_{i\neq k} \vert \tilde{b}_{ii} \vert^2 \right).
\end{split}
\end{equation*}

Similarly we obtain
\begin{equation*}
(\mathrm{III})  \leq  \frac{1}{n^2} \left(\sum_{i \neq j} ( \vert \tilde{c}_{ij}\vert^2 +\gamma^2 \bar{\tilde{c}}_{ij} \tilde{c}_{ji}) + (\kappa^2(n)- \gamma^2) \sum_{i \neq k} \vert \tilde{c}_{ii} \vert^2 \right), \end{equation*}
and 
\begin{equation*}
(\mathrm{IV})  \leq \frac{1}{n^2} \left( \sum_{i \neq j} (\vert \tilde{a}_{ij} \vert^2 + \bar{\tilde{a}}_{ij} \tilde{a}_{ji}  ) + (\kappa^2(n)-1) \sum_{i \neq k} \vert \tilde{a}_{ii}\vert^2 \right).
\end{equation*}

Note that 
\begin{equation*}
\sum_{i \neq j} (\vert \tilde{a}_{ij} \vert^2 +  \vert \tilde{b}_{ij} \vert^2 + \vert \tilde{c}_{ij} \vert^2 + \vert \tilde{d}_{ij} \vert^2 ) \leq \sum_{i, j \neq k} {\rm Tr} \left[ (\widetilde{R}_{ij}^{(k)})^* \widetilde{R}_{ij}^{(k)} \right],
\end{equation*}
\begin{equation*}
\sum_{i \neq j} (\bar{\tilde{a}}_{ij} \tilde{a}_{ji} + \gamma^2 \bar{\tilde{b}}_{ij} \tilde{b}_{ji} + \gamma^2 \bar{\tilde{b}}_{ij} \tilde{b}_{ji} +  \bar{\tilde{d}}_{ij} \tilde{d}_{ji} ) \leq \sum_{i, j \neq k}{\rm Tr} \left[ \Sigma( \widetilde{R}_{ij}^{(k)})^* \Sigma(\widetilde{R}_{ji}^{(k)}) \right],
\end{equation*}
and
\begin{equation*}
 \sum_{i \neq k} (\kappa^2(n)-1) ( \vert \tilde{a}_{ii}\vert^2 + \vert \tilde{d}_{ii}\vert^2 ) + (\kappa^2(n)- \gamma^2)  (\vert \tilde{c}_{ii} \vert^2 +  \vert \tilde{b}_{ii} \vert^2) \leq \kappa^2(n)\sum_{i \neq k}  {\rm Tr} \left[ (\widetilde{R}_{ii}^{(k)})^* \widetilde{R}_{ii}^{(k)} \right].
\end{equation*}

In summary, we have 
\begin{equation*}
\begin{split}
{\rm Tr} (Z_k) 
& \leq \frac{1}{n^2}  {\rm Tr}\left[ \sum_{i,j \neq k}( (\widetilde{R}_{ij}^{(k)})^* \widetilde{R}_{ij}^{(k)} + \Sigma( \widetilde{R}_{ij}^{(k)})^* \Sigma(\widetilde{R}_{ji}^{(k)}))  + \kappa^2(n)\sum_{i \neq k} (\widetilde{R}_{ii}^{(k)})^* \widetilde{R}_{ii}^{(k)} \right]\\
& = \frac{1}{n^2}  \left( {\rm Tr} ((\widetilde{R}^{(k)})^* \widetilde{R}^{(k)}) + {\rm Tr} (\Sigma(\widetilde{R}^{(k)})^* \Sigma(\widetilde{R}^{(k)})) +  \kappa^2(n)\sum_{i \neq k} (\widetilde{R}_{ii}^{(k)})^*\widetilde{R}_{ii}^{(k)}  \right)\\
& \leq \frac{2}{n} (\Vert \widetilde{R}^{(k)} \Vert^2 + \Vert \Sigma (\widetilde{R}^{(k)})\Vert^2) + \frac{1}{n} \kappa^2(n) \Vert \widetilde{R}^{(k)} \Vert^2\\
& \leq \frac{2(1+C^2)+ \kappa^2(n)}{n}  \frac{1}{\Im (\eta)^2},
\end{split}
\end{equation*}
where we $C>0$ is a constant such that $\Vert \Sigma \Vert \leq C$. It follows that $\|Z_k\|  \leq {\rm Tr}(Z_k) \rightarrow 0$ as $n \rightarrow \infty.$
\end{proof}

\begin{proof}[Proof of Theorem \ref{thm:fix-point-equ}]
Recall that for any $k=1, \ldots, n,$ by the Schur block inversion formula, we have 
\begin{equation}\label{eq:inversion}
R(q)_{kk} = \left( \begin{bmatrix} 0 &
M_{kk} \\  \bar{M}_{kk} & 0 \end{bmatrix} -q - {Q^{(k)}}^* R^{(k)}(q) Q^{(k)} \right)^{-1}.
\end{equation}

We first claim that the norm $\Vert {Q^{(k)}}^* Q^{(k)} \Vert$ centers around its expectation in probability. We have
\begin{equation*}
 {Q^{(k)}}^* Q^{(k)} = \frac{1}{n}  \begin{bmatrix} \sum_{i \neq k} \vert X_{ik} \vert^2 &
0 \\  0 & \sum_{i \neq k} \vert X_{ki} \vert^2 \end{bmatrix}.
\end{equation*}
Thus, it suffices to show that $\frac{1}{n} \sum_{i \neq k} \vert X_{ik} \vert^2$ centers around $\mathbb{E} [\frac{1}{n} \sum_{i \neq k} \vert X_{ik} \vert^2].$ To this end, 
we consider
\begin{equation*}
\begin{split}
\mathbb{E} \vert \frac{1}{n}\sum_{i \neq k} \vert X_{ik} \vert^2 - \mathbb{E} [\frac{1}{n} \sum_{i \neq k} \vert X_{ik} \vert^2] \vert^2 & = \frac{1}{n^2} \mathbb{E} [\sum_{i,j \neq k} (
\vert X_{ik} \vert^2  \vert X_{jk} \vert^2 - \mathbb{E} \vert X_{ik} \vert^2   \cdot \mathbb{E} \vert X_{jk} \vert^2) ] \\
&  \leq \frac{1}{n^2} \mathbb{E}[\sum_{i \neq k} \vert X_{ik} \vert^4]\\
& \leq \frac{\kappa^2(n)}{n} \rightarrow 0, \; \text{as} \; n \rightarrow \infty.
\end{split}
\end{equation*}
Our claim follows from the Markov inequality. 

Let $\varepsilon_1 = {Q^{(k)}}^* R^{(k)} Q^{(k)} - {Q^{(k)}}^* \widetilde{R}^{(k)} Q^{(k)},$ then by the resolvent identity we have 
\begin{equation*}
\begin{split}
\Vert \varepsilon_1 \Vert &  \leq  \Vert {Q^{(k)}}^* Q^{(k)} \Vert \cdot  \Vert R^{(k)} - \widetilde{R}^{(k)} \Vert  \\
&=  \Vert {Q^{(k)}}^* Q^{(k)} \Vert \cdot \Vert \widetilde{R}^{(k)} (B^{(k)} - \widetilde{B}^{(k)}) R^{(k)} \Vert  \\
& \leq \frac{ \Vert {Q^{(k)}}^* Q^{(k)} \Vert } {\sqrt{n} \Im (\eta)^2} \cdot \left\| {\rm diag} \; \left(  \begin{bmatrix} 0 &
X_{ki} \\  \bar{X}_{ik} & 0 \end{bmatrix} \right)_{i \neq k}\right\|\\
&  \leq \frac{\kappa(n) \Vert {Q^{(k)}}^* Q^{(k)}\Vert}{\sqrt{n} \Im (\eta)^2}.
\end{split}
\end{equation*}
Recall that  $\Vert {Q^{(k)}}^* Q^{(k)} \Vert \rightarrow 1$ in probability as $n \rightarrow \infty$. Consequently, $\Vert \varepsilon_1 \Vert $ converges to $0$ in probability as $n \rightarrow \infty.$ Thus, we obtain
\begin{equation}\label{eq:inversion-R'}
R_{kk}  = \left( \begin{bmatrix} 0 &
M_{kk} \\  \bar{M}_{kk} & 0 \end{bmatrix} - q - {Q^{(k)}}^* \widetilde{R}^{(k)} Q^{(k)} - \varepsilon_1 \right)^{-1},
\end{equation}
where $ \Vert \varepsilon_1 \Vert = o(1)$ in probability. 

Combining Lemma \ref{lem:epsilon2}, Lemma \ref{lem:epsilon3}, and Lemma \ref{lem:epsilon4} {\rm (1)}, we have 
\begin{equation*}
\begin{split}
&\left\Vert \mathbb{E}_{k} \left[ {Q^{(k)}}^* \widetilde{R}^{(k)} Q^{(k)} \right]  - \Sigma(\Gamma_n(q)) \right\Vert  = \left \Vert \Sigma (\widetilde{\Gamma}_n^{(k)}(q))- \Sigma(\Gamma_n(q)) \right\Vert\\
& \leq \left\Vert  \Sigma (\widetilde{\Gamma}_n^{(k)} (q)) - \Sigma (\Gamma^{(k)}_n(q)) \right\Vert + \left\Vert \Sigma (\Gamma^{(k)}_n(q)) - \Sigma (\Gamma_n(q)) \right\Vert\\
& \leq C \cdot \left\Vert  \widetilde{\Gamma}_n^{(k)} (q)) - \Gamma^{(k)}_n(q) \right\Vert + C \cdot \left\Vert \Gamma^{(k)}_n(q) - \Gamma_n(q) \right\Vert\\
& = O\left(\frac{\kappa(n)}{ \sqrt{n} \Im (\eta)^{2} } \right) + O\left(\frac{1}{ n \Im (\eta) }\right)\rightarrow 0 \; \text{in probability as} \; n \rightarrow \infty,
\end{split}
\end{equation*}
where $C>0$ is a constant such that $\Vert \Sigma \Vert \leq C$. Moreover, it follows from Lemma \ref{lem:epsilon4} {\rm (2)} and Lemma \ref{lem:epsilon5} that 
\begin{equation*}
\begin{split}
& \left\Vert {Q^{(k)}}^* \widetilde{R}^{(k)} Q^{(k)} - \mathbb{E} \left[ \Sigma(\Gamma_n(q)) \right] \right\Vert \\
& \leq \left\Vert {Q^{(k)}}^* \widetilde{R}^{(k)} Q^{(k)} -  \Sigma(\Gamma_n(q)) \right\Vert + \left\Vert \Sigma(\Gamma_n(q)) - \mathbb{E} \left[ \Sigma(\Gamma_n(q)) \right] \right\Vert\\
& \rightarrow 0 \; \text{in probability as} \; n \rightarrow \infty.
\end{split}
\end{equation*}

We denote 
\begin{equation*}
A = -  \begin{bmatrix} 0 &
M_{kk} \\  \bar{M}_{kk} & 0 \end{bmatrix} + \begin{bmatrix} 0 &
G \\  \overline{G} & 0 \end{bmatrix} + {Q^{(k)}}^* \widetilde{R}^{(k)} Q^{(k)}  - \mathbb{E} \left[ \Sigma(\Gamma_n(q)) \right].
\end{equation*}
Using the Lindeberg--Feller central limit theorem and the Skorokhod representation theorem \cite{Resnick}, we have almost surely as $n \rightarrow \infty$ 
\begin{equation*}
\left\|  \begin{bmatrix} 0 &
M_{kk} \\  \bar{M}_{kk} & 0 \end{bmatrix} - \begin{bmatrix} 0 &
G \\  \overline{G} & 0 \end{bmatrix}  \right\|  \rightarrow 0. 
\end{equation*}
Therefore, the $2 \times 2$ matrix $A$ has a norm that converges to $0$ in probability as $n \rightarrow \infty.$

By \eqref{eq:inversion-R'}, 
\begin{equation*}
R_{kk} A = -R_{kk} \left( q +  \mathbb{E} \left[ \Sigma(\Gamma_n(q)) \right] -  \begin{bmatrix} 0 &
G \\  \overline{G} & 0 \end{bmatrix} +\varepsilon_1 \right) - \un_2,
\end{equation*}
which yields that 
\begin{equation*}
R_{kk} = - \left( A+ q +  \mathbb{E} \left[ \Sigma(\Gamma_n(q)) \right]  -  \begin{bmatrix} 0 &
G \\  \overline{G} & 0 \end{bmatrix} + \varepsilon_1 \right)^{-1}
\end{equation*}
with $\| \varepsilon_1\| = o(1)$ in probability.

For any $k$, denote a sequence of random variables $\{Z_{k,n}\}_{n \geq 1}$ by 
\begin{equation*}
Z_{k, n} :=   \left\| R_{kk} - \left(  -q -  \mathbb{E} \left[ \Sigma(\Gamma_n(q)) \right] +  \begin{bmatrix} 0 &
G \\  \overline{G} & 0 \end{bmatrix} \right)^{-1}  \right\|. 
\end{equation*}
Since the norms of  $$\left( q +  \mathbb{E} \left[ \Sigma(\Gamma_n(q)) \right]-  \begin{bmatrix} 0 &
G \\  \overline{G} & 0 \end{bmatrix}  \right)^{-1}$$
and $R_{kk}$ are bounded above by  $1/ \Im (\eta),$ $\{ Z_{k,n} \}_{n \geq 1}$ is uniformly integrable. Moreover, by the resolvent identity, 
 \begin{equation*}
\begin{split}
Z_{k, n} & =  \left\| R_{kk}  \cdot (A + \varepsilon_1) \cdot \left( -q -  \mathbb{E} \left[ \Sigma(\Gamma_n(q)) \right] +  \begin{bmatrix} 0 &
G \\  \overline{G} & 0 \end{bmatrix} \right)^{-1}   \right\| \\
& \leq \frac{1}{\Im (\eta)^2} \cdot ( \|A\|+ \|\varepsilon_1 \|). 
\end{split}
\end{equation*}
Thus, $Z_{k,n} \rightarrow 0$ in probability as $n \rightarrow \infty.$ Combining the uniform integrability, we have $\mathbb{E} [Z_{k,n}] \rightarrow 0$ as $n \rightarrow \infty.$
Therefore we obtain
\begin{equation*}
\mathbb{E} R_{kk} = - \mathbb{E}  \left( q +  \mathbb{E} \left[ \Sigma(\Gamma_n(q)) \right]-  \begin{bmatrix} 0 &
G \\  \overline{G} & 0 \end{bmatrix}  \right)^{-1} + \varepsilon
\end{equation*}
with $\|\epsilon\| = o(1).$ Note that $\Gamma_n(q)  = \frac{1}{n} \sum_{k=1}^n R_{kk}.$  Thus, 
\begin{equation}
\mathbb{E} \left[ \Gamma_n(q) \right] = - \mathbb{E}  \left( q + \mathbb{E} \left[ \Sigma(\Gamma_n(q)) \right]-  \begin{bmatrix} 0 &
G \\  \overline{G} & 0 \end{bmatrix}  \right)^{-1} + \varepsilon.
\end{equation}
This finishes the proof. 
\end{proof}

\subsection{Uniqueness of the fixed point equation and some free probability results}
In this subsection, we show that the fixed point equation \eqref{eq:fixed-point-eq} obtained in Section \ref{sec:approx-fixed-pt} has a unique solution and such a solution has a natural interpretation in free probability. 

\subsubsection{Operator-valued free probability}
Let $a,  y \in \mathcal{A}$ such that $a$ is $*$-freely independent from $y$ in $(\mathcal{A}, \phi).$ We write $A=\begin{bmatrix}
0 & a\\
a^* & 0
\end{bmatrix}$ and $B=\begin{bmatrix}
0 & y\\
y^* & 0
\end{bmatrix}$. Then $A, B$ live in the operator-valued probability space $(\mathcal{M}_2(\mathcal{A}), \Phi, \mathcal{M}_2(\mathbb{C}))$, where $\mathcal{M}_2(\mathcal{A})$ is the $2\times 2$ matrices with entries from $\mathcal{A}$ and $\Phi := \phi\otimes \un_2: \mathcal{M}_2(\mathcal{A})\rightarrow \mathcal{M}_2(\mathbb{C})$ is the entrywise conditional expectation. It is known that the matrix-valued operators $A, B$ are free with amalgamation over $\mathcal{M}_2(\mathbb{C})$ in the operator-valued free probability space $(\mathcal{M}_2(\mathcal{A}), \Phi, \mathcal{M}_2(\mathbb{C}))$ (see \cite[Chapter 9]{MingoSpeicherBook} for example). 

If $a\in \log^+(\mathcal{A})\subset \widetilde{\mathcal{A}}$, then $A$ is affiliated with $\mathcal{M}_2({\mathcal{A}})$. 
For any $q=q(z,\eta)\in \mathbb{H}_+$, the operator-valued Cauchy transform of $A$ is defined as
\begin{equation}
G_A(q)=\Phi \left( (q-A)^{-1}\right),
\end{equation}
and the operator Stieltjes transform of $A$ is defined as
\begin{equation}
S_A(q)=-G_A(q)=\Phi \left( (A-q)^{-1}\right).
\end{equation}
If $a\in\mathcal{A}$ is a bounded operator, then the $R$-transform of $A$ is defined as
\begin{equation}
R_A(b)=G^{\langle-1\rangle}(b)-b^{-1},
\end{equation}
provided that $b$ is invertible and $\Vert  b \Vert$ is sufficiently small. 
There are a pair of two analytic functions
$\Omega_1, \Omega_2: \mathbb{H}_+\rightarrow \mathbb{H}_+$, called subordination functions, such that
\[
G_{A+B}(q)=G_A(\Omega_1(q))=G_B(\Omega_2(q)), \; q\in\mathbb{H}_+.
\]
Moreover, for any $q\in\mathbb{H}_+$, we have 
\[
\Omega_1(q)+\Omega_2(q)=q+(G_{A+B}(q))^{-1}.
\]
The reader is referred to \cite{BelinschiTR2018-sub-operator-valued, SpeicherAMS1998, DVV-operator-valued-1992} for a complete treatment. Since the operator $a$ is unbounded, we have to use the subordination functions for unbounded operators that was established in \cite[Section 3]{BelinschiYinZhong2021Brown} (under some assumptions).

\subsubsection{Operator-valued Stieltjes transform of the sum}

In what follows, we denote $A=\begin{bmatrix}
0 & a\\
a^* & 0
\end{bmatrix}$ and $B=\begin{bmatrix}
0 & y\\
y^* & 0
\end{bmatrix}$, where $a$ is a Gaussian distributed normal operator with law $\mathcal{N}(0, K_\mu)$ and $y \in \{c, g_{\gamma}\}$ freely independent from $a$.

The following result was essentially proved in \cite[Theorem 3.8]{Zhong2021Brown_ctgt} for bounded case and \cite[Theorem 7.3]{BelinschiYinZhong2021Brown} for unbounded case. Since the proofs uses the operator-valued subordination functions for unbounded operators \cite[Theorem 3.1]{BelinschiYinZhong2021Brown} and we only need a special form of the general result \cite[Lemma 4.3]{BelinschiYinZhong2021Brown},  we include a short proof for reader's convenience. 

\begin{theorem}[The fixed point equation, \cite{Zhong2021Brown_ctgt, BelinschiYinZhong2021Brown}]
	\label{thm:subordination-main}
	For $q=q(z,\eta)\in \mathbb{H}_+$, the operator-valued Stieltjes transform 
	\[
	S_{A+B}(q)=- \Phi \big( (q-A-B)^{-1}\big)
	\]
	is a solution of the fixed point equation
	\begin{equation}\label{eqn:fixed-pt-eqn-Cauchy}
	\begin{bmatrix}
	\alpha  & \beta\\
	\bar{\beta} &\alpha
	\end{bmatrix}
	=- \Phi \left( q + \Sigma_{\gamma} \begin{bmatrix}
	\alpha  & \beta\\
	\bar{\beta} &\alpha
	\end{bmatrix}    - A  \right) ^{-1},
	\end{equation}
	where 
	\[
	\Sigma_{\gamma} \begin{bmatrix}
	\alpha  & \beta\\
	\bar{\beta} &\alpha
	\end{bmatrix} 
	=\begin{bmatrix}
	 \alpha  &  \gamma \bar{\beta}\\
	\bar{\gamma} \beta &  \alpha 
	\end{bmatrix}.
	\]
\end{theorem}

\begin{proof}
	Since $A, B$ are free with amalgamation in the operator-valued probability space $(\mathcal{M}_2(\mathcal{A}), \Phi, \mathcal{M}_2(\mathbb{C}))$, we have 
	\[
	\Omega_1(q)+\Omega_2(q)=q+(G_{A+B}(q))^{-1}.
	\]
	From the subordination relation $G_{A+B}(q)=G_B(\Omega_2(q))$ with $q=q(z,\eta)$, if $\Im z$ is sufficiently large, then $\Vert G_{A+B}(q) \Vert\leq 1/ \Im z$ is small, and hence for such $q$ we have 
	\begin{equation*}
	\Omega_2(q)=R_B (G_{A+B}(q))+(G_{A+B}(q))^{-1}.
	\end{equation*}
	Consequently, we have
	\begin{equation}\label{eqn:formula-Omega1}
	\Omega_1(q)=q-R_B (G_{A+B}(q))=q+R_B (S_{A+B}(q)).
	\end{equation}
	It is known that the operator $B$ is an operator-valued semicircular element (see \cite[Proposition 2.1]{Zhong2021Brown_ctgt} or \cite[Section 9.5]{MingoSpeicherBook})) whose $R$-transform can be explicitly written as
	\begin{equation}\label{eqn:formula-R-transform-Y}
	R_B \left( \begin{bmatrix}
	\alpha  & \beta\\
	\bar{\beta}  &\alpha
	\end{bmatrix}   \right) =\begin{bmatrix}
	\kappa (y, y^*) \alpha &  \kappa(y,y) \bar{\beta}\\
	\kappa(y^*, y^*) \beta  &  \kappa(y^*, y) \alpha
	\end{bmatrix}=
	\begin{bmatrix}
	 \alpha  & \bar{\beta} \gamma \\
	\beta  \bar{\gamma}    &  \alpha 
	\end{bmatrix},
	\end{equation}
	where we used the free cumulant formula \eqref{eqn:free-cumulnts-y}. 
	
	Note that, for any $q\in \mathbb{H}_+$, we have $S_{A+B}(q)\in \mathbb{H}_+$. Hence, for any $q\in\mathbb{H}_+$, the right hand side of \eqref{eqn:formula-Omega1} is in $\mathbb{H}_+$. By the uniqueness of the subordination function (see \cite{BelinschiTR2018-sub-operator-valued} and \cite[Theorem 3.1]{BelinschiYinZhong2021Brown}) and its analytic extension, we conclude that the formula for the subordination function $\eqref{eqn:formula-Omega1}$ holds for any $q\in\mathbb{H}_+$. Hence, we can rewrite it as
	\[
	\Omega_1(q)=q-\Sigma_{\gamma} (G_{A+B}(q))=q+\Sigma_{\gamma} (S_{A+B}(q)).
	\]
	Therefore, for any $q\in\mathbb{H}_+$, we obtain 
	\[
	S_{A+B}(q)=-G_{A+B}(q)=-G_A(\Omega(q))=-G_A(q+\Sigma_{\gamma}(S_{A+B}(q))), 
	\]
	which exactly means that the Stieltjes transform $S_{A+B}(q)$ satisfies the fixed point equation \eqref{eqn:fixed-pt-eqn-Cauchy}.
\end{proof}

\subsubsection{The uniqueness of the fixed point equation}
We next show the fixed point equation \eqref{eqn:fixed-pt-eqn-Cauchy} has a unique solution.
For any $\varepsilon>0$, denote a subset $\mathbb{H}_+^\varepsilon$ of $\mathbb{H}_+$ by 
\[
\mathbb{H}_+^\varepsilon=\left\{ \begin{bmatrix}
\sqrt{-1}\varepsilon & z\\
\overline{z} & \sqrt{-1}\varepsilon
\end{bmatrix}: z\in\mathbb{C}    \right\}.
\]

\begin{lemma}\label{lemma:unique}
	For any $\varepsilon>0$, 
	the fixed point equation \eqref{eqn:fixed-pt-eqn-Cauchy} has a unique solution in $\mathbb{H}_+^\varepsilon$ for any $q=q(z,\sqrt{-1}\varepsilon)\in\mathbb{H}_+$ with $\alpha=\alpha(q)\in \sqrt{-1}\mathbb{R}_+$. The solution is given by
	\begin{equation}\nonumber
		\begin{bmatrix}
		\alpha & \beta\\
		\bar{\beta}  & \alpha
		\end{bmatrix}
		  =S_{A+B}(q)=- \Phi \big( (q-A-B)^{-1}\big). 
	\end{equation}
In particular, 
\begin{equation}\label{eqn:alpha-Cauchy}
	 \alpha=s_{A+B, 11}(z,\varepsilon)=\sqrt{-1} \varepsilon \phi\bigg(\big( (z-a-y)(z-a-y)^*+\varepsilon^2\big)^{-1} \bigg).
\end{equation}
\end{lemma}
\begin{proof}
	Let $\begin{bmatrix}
	\alpha & \beta\\
	\bar{\beta}  & \alpha
	\end{bmatrix}$ be a solution for the fixed point equation \eqref{eqn:fixed-pt-eqn-Cauchy} for $q=q(z,\sqrt{-1}\varepsilon)$. 
	Define the function $\Omega:\mathbb{H}_+^\varepsilon\rightarrow \mathbb{H}_+^\varepsilon$ by
	\[
	\Omega\left( q \right)=q+\Sigma_\gamma \begin{bmatrix}
	\alpha & \beta\\
	\bar{\beta} & \alpha
	\end{bmatrix}, \;
	q=q(z, \sqrt{-1}\varepsilon).
	\]
	That is, 
	\[
	\Omega\left( \begin{bmatrix}
	\sqrt{-1}\varepsilon & z\\
	\bar{z} & \sqrt{-1}\varepsilon
	\end{bmatrix}   \right)=\begin{bmatrix}
	\sqrt{-1}\varepsilon + \alpha & z+\gamma \bar{\beta} \\
	\bar{z}+ \bar{\gamma} \beta & \sqrt{-1}\varepsilon+\alpha
	\end{bmatrix}.
	\]
	Then, by Equation \eqref{eqn:fixed-pt-eqn-Cauchy}, $\Omega$ can be rewritten as
	\begin{equation} \label{eqn:Omega}
	\Omega(q)=q-\Sigma_\gamma ( G_A (\Omega(q)) ).
	\end{equation}
	Denote by 
	\[
	H(b):=b+\Sigma_\gamma(G_A(b)).
	\]
	Then, Equation \eqref{eqn:Omega} can be reorganized as
	\[
	H(\Omega(q))=q.
	\]

	We will study the image $\Omega(\mathbb{H}_+^\varepsilon)$ and properties of its left inverse $H$. Let $\lambda=z+\gamma \bar{\beta}$. By comparing the $(1,1)$-entry of both sides of \eqref{eqn:fixed-pt-eqn-Cauchy}, for $q\in \mathbb{H}_+^\varepsilon$, 
	\[
	\alpha=(\sqrt{-1}\varepsilon+\alpha)\phi \big( ((\lambda-a)^*(\lambda-a)-(\sqrt{-1}\varepsilon+\alpha)^2 )^{-1} \big).
	\]
	Let $w=\Im (\alpha(q))+\varepsilon$, we obtain
	\[
	w-\varepsilon=w\phi\big( ((\lambda-a)^*(\lambda-a)+w^2 )^{-1} \big).
	\]
	Hence, $w>\varepsilon$ and 
	\begin{equation} \label{eqn:fixed-pt-h}
	1=\frac{w}{w-\varepsilon} \phi\big( ((\lambda-a)^*(\lambda-a)+w^2 )^{-1} \big).
	\end{equation}
	The right-hand side of the above equation is a decreasing function of $w$ on $(\varepsilon,\infty)$ with range $(0,\infty)$. Thus, there is a unique $w>\varepsilon$ depending only on $\lambda$ and $\varepsilon$ satisfying \eqref{eqn:fixed-pt-h}. We denote this unique
	solution for \eqref{eqn:fixed-pt-h} by $w^{(\lambda)}(\varepsilon)$. Hence, 
	\[
	\Omega\left( \begin{bmatrix}
	\sqrt{-1}\varepsilon & z\\
	\bar{z} & \sqrt{-1}\varepsilon
	\end{bmatrix}   \right)=\begin{bmatrix}
	\sqrt{-1} w^{(\lambda)}(\varepsilon) & \lambda\\
	\bar{\lambda}  &  \sqrt{-1} w^{(\lambda)}(\varepsilon)
	\end{bmatrix},
	\]
	where $\lambda=z+\gamma\bar{\beta}$, 
	and the fixed point equation \eqref{eqn:fixed-pt-eqn-Cauchy} is rewritten as 
	\begin{equation}
	\begin{bmatrix}
	\alpha  & \beta\\
	\bar{\beta} &\alpha
	\end{bmatrix}
	=- \Phi \left( \begin{bmatrix}
	\sqrt{-1}  w^{(\lambda)}(\varepsilon) & \lambda\\
	\bar{\lambda}  & \sqrt{-1} w^{(\lambda)}(\varepsilon)
	\end{bmatrix}  - A \right) ^{-1}.
	\end{equation}
	By comparing $(1,1)$ and $(1,2)$-entries, we have 
	\begin{equation}\label{eqn:subordination-general}
	\begin{aligned}
	\alpha&=\sqrt{-1}  w^{(\lambda)}(\varepsilon)\phi\big( ((\lambda-a)^*(\lambda-a)+(w^{(\lambda)}(\varepsilon))^2 )^{-1} \big),\\
	\beta&=-\phi \bigg( (\lambda-a)\big( (\lambda-a)^*(\lambda-a)+(w^{(\lambda)}(\varepsilon))^2  \big)^{-1}  \bigg).
	\end{aligned}
	\end{equation}

	In summary, if $\lambda=z+\gamma \bar{\beta}$, then $\sqrt{-1}\varepsilon+\alpha=\sqrt{-1} w^{(\lambda)}(\varepsilon).$ 
	In other words, 
	\begin{equation} \label{eqn:image-Omega}
	\Omega (\mathbb{H}_+^\varepsilon)\subset \left\{ \begin{bmatrix}
	\sqrt{-1}  w & \lambda\\
	\bar{\lambda} & \sqrt{-1}  w
	\end{bmatrix} : w= w^{(\lambda)}(\varepsilon) , \lambda\in\mathbb{C}  \right\}.
	\end{equation}
	Fix $\varepsilon>0$, we hence focus on the relation between $z$ and $\lambda=z+\gamma\bar{\beta}$. 
	
	We now turn to the study of the solution for the fixed point equation \eqref{eqn:fixed-pt-eqn-Cauchy}  provided by the operator-valued Stieltjes transform 
	$S_{A+B}(q)$ from Theorem \ref{thm:subordination-main}. This was known in \cite[Theorem 3.8 and Proposition 5.1]{Zhong2021Brown_ctgt}. We include a somewhat direct argument for convenience. We have 
	\begin{align*}
	S_{A+B}(q)&=\Phi (A+B-q)^{-1}= \Phi \begin{bmatrix}
	- \sqrt{-1} \varepsilon & a+y-z\\
	(a+y-z)^* & -\sqrt{-1} \varepsilon
	\end{bmatrix}^{-1}\\
	&=\begin{bmatrix}
	s_{A+B, 11}(z,\varepsilon) &  s_{A+B, 12}(z,\varepsilon) \\
	s_{A+B, 21}(z,\varepsilon)  &  s_{A+B, 22}(z,\varepsilon) 
	\end{bmatrix},
	\end{align*}
	where 
	\begin{align*}
	s_{A+B, 11}(z,\varepsilon)&=\sqrt{-1} \varepsilon \phi\bigg(\big( (z-a-y)(z-a-y)^*+\varepsilon^2\big)^{-1} \bigg),\\
	s_{A+B, 22}(z,\varepsilon)&=\sqrt{-1} \varepsilon \phi\bigg(\big( (z-a-y)^*(z-a-y)+\varepsilon^2\big)^{-1} \bigg),\\
	s_{A+B, 12}(z,\varepsilon)&=-\phi\bigg((z-a-y)\big( (z-a-y)^*(z-a-y)+\varepsilon^2\big)^{-1} \bigg),\\
	s_{A+B, 21}(z,\varepsilon)&=-\phi\bigg((z-a-y)^*\big( (z-a-y)(z-a-y)^*+\varepsilon^2\big)^{-1} \bigg).
	\end{align*}
	By the tracial property, we have $s_{A+B, 11}(z,\varepsilon)=s_{A+B, 22}(z,\varepsilon)$ and 
	$s_{A+B, 21}(z,\varepsilon)=\overline{s_{A+B, 12}(z,\varepsilon)}$. 
	Hence, we also have 
	\[
	\Omega_1\left( \begin{bmatrix}
	\sqrt{-1}\varepsilon & z\\
	\overline{z} & \sqrt{-1}\varepsilon
	\end{bmatrix}   \right)=\begin{bmatrix}
	\sqrt{-1}  w^{(\lambda)}(\varepsilon) & \lambda\\
	\overline{\lambda}  & \sqrt{-1}  w^{(\lambda)}(\varepsilon)
	\end{bmatrix},
	\]
	where 
	\begin{equation}
	\lambda=z-\gamma \phi\bigg((z-a-y)^*\big( (z-a-y)(z-a-y)^*+\varepsilon^2\big)^{-1} \bigg).
	\end{equation}
	Equation \eqref{eqn:subordination-general} is exactly the entries of the subordination relation
	$-S_{A+B}(q)=-S_A(\Omega_1(q))$. More precisely,
	\begin{align*}\label{eqn:subordination-entries}
	\sqrt{-1}\varepsilon \phi&\bigg(\big( (z-a-y)(z-a-y)^*+\varepsilon^2\big)^{-1} \bigg)\\
	&\qquad\qquad 
	=  \sqrt{-1} w^{(\lambda)}(\varepsilon)\phi\big( ((\lambda-a)^*(\lambda-a)+(w^{(\lambda)}(\varepsilon))^2 )^{-1} \big),\\
	-\phi&\bigg((z-a-y)\big( (z-a-y)^*(z-a-y)+\varepsilon^2\big)^{-1} \bigg)\\
	&\qquad\qquad=-\phi \bigg( (\lambda-a)\big( (\lambda-a)^*(\lambda-a)+(w^{(\lambda)}(\varepsilon))^2  \big)^{-1}  \bigg).
	\end{align*}

	We now define a continuous function $J^{(\varepsilon)}:\mathbb{C}\rightarrow\mathbb{C}$ by 
	\[
	J^{(\varepsilon)}(z)=z-\gamma \phi\bigg((z-a-y)^*\big( (z-a-y)(z-a-y)^*+\varepsilon^2\big)^{-1} \bigg).
	\]
	It is clear that
	\[
	\lvert J^{(\varepsilon)}(z)-z \rvert \leq \vert \gamma\vert/\varepsilon.
	\]
	In particular $J^{(\varepsilon)}(\infty)=\infty$. On the other hand, $J^{(\varepsilon)}$ is an injective map; see \cite[Lemma 4.10]{BelinschiYinZhong2021Brown}. Hence $J^{(\varepsilon)}$ can be considered as a $C^\infty$ map from $\mathbb{C} \cup \{\infty\} = \mathbb{S}^2$ to itself. By applying  Borsuk-Ulam Theorem \cite{RotmanBook} for dimension two, we can show that $ J^{(\varepsilon)}$ is surjective. Therefore, 
	\[
	\Omega_1 (\mathbb{H}_+^\varepsilon)= \left\{ \begin{bmatrix}
	\sqrt{-1} w & \lambda\\
	\bar{\lambda} &  \sqrt{-1}w
	\end{bmatrix} : w= w^{(\lambda)}(\varepsilon) , \lambda\in\mathbb{C}  \right\}.
	\]
	Finally, we finish the proof by showing that $\Omega(q)=\Omega_1(q)$ for any $q=q(z,\sqrt{-1}\varepsilon)\in\mathbb{H}_+^\varepsilon$, which would yield that 
	\[
	\alpha=s_{A+B, 11}(z,\varepsilon), \; \beta=s_{A+B, 12}(z,\varepsilon).
	\]
	Indeed, for any $b$ of the form $\begin{bmatrix}
	\sqrt{-1} w^{(\lambda)}(\varepsilon) & \lambda\\
	\bar{\lambda}  & \sqrt{-1} w^{(\lambda)}(\varepsilon)
	\end{bmatrix}$ with $\lambda\in\mathbb{C}$ and $w^{(\lambda)}(\varepsilon)$ defined as the unique solution for \eqref{eqn:fixed-pt-h}, we have 
	\[
	\Omega_1(H(b))=b.
	\]
	Therefore, by plugging $b=\Omega(q)$, we have 
	\[
	\Omega(q)=\Omega_1(H(\Omega(q)))=\Omega_1(q),
	\]
	where we used $H(\Omega(q))=q$ for any $q=q(z,\sqrt{-1}\varepsilon)$. 
\end{proof}

\begin{remark}\label{remark:connection-subordination-w}
	The function $w^{(\lambda)}(\cdot)$ is related to the scalar-valued subordination function. Given any $\lambda\in\mathbb{C}$, denote by $\mu_1=\widetilde{\mu}_{\vert x-\lambda \vert}$ and $\mu_2$ be the standard semicircular distribution. Let $\omega_1$ be the subordination function such that
	\[
	G_{\mu_1\boxplus\mu_2}(z)=G_{\mu_1}(\omega_1(z)), \qquad z\in\mathbb{C}^+.
	\]
	Then, $w^{(\lambda)}(\varepsilon)=\Im \omega_1(i\varepsilon)$. See \cite[Proposition 3.5]{Zhong2021Brown_ctgt}.
\end{remark}

\subsection{Completing the proof of Theorem \ref{thm:singular-value}}
Recall that $H(z)$ is a $2n \times 2n$ Hermitian matrix given as follows:
\begin{equation*}
H(z) := \begin{bmatrix} 0 &
M-z \\  (M-z)^* & 0 \end{bmatrix}.
\end{equation*}
For $\eta \in \mathbb{C}_+$, from \eqref{defn:S-transform-singular-V2}, the Stieltjes transform of $\widetilde{\nu}_{M-z}$ can be calculated as  
\begin{equation}
S_{\widetilde{\nu}_{M-z}} (\eta)= \frac{1}{2n} {\rm Tr} \left[ (H(z)- \eta \un_{2n})^{-1} \right]=\frac{1}{2n} {\rm Tr}\left[ R(q) \right]=a(q).
\end{equation}

\begin{lemma}[Tightness]\label{lem:tight}
The averaged measure $\mathbb{E} [\nu_{M -z}]$ is tight. 
\end{lemma}
\begin{proof}
It suffices to show that 
\begin{equation}
\mathbb{E} \int s^2 d\nu_{M-z}(s) < \infty. 
\end{equation}
By the Weyl inequality, we have
\begin{equation*}
\begin{split}
\mathbb{E} \int s^2 d\nu_{M-z}(s)  & \leq  \mathbb{E} \int s^2 d\nu_M(s) + |z|^2 =  \mathbb{E} {\rm tr} [\vert M \vert^2] +  \vert z \vert^2 \\
& = \frac{1}{n} \sum_{i, j =1}^n \mathbb{E} [ \vert M_{ij} \vert^2] + \vert z \vert^2 \\
& = \frac{1}{n^2} \sum_{i \neq k} \mathbb{E} [\vert X_{ik} \vert^2 ] + \frac{1}{n^2} \sum_{i=1}^n \mathbb{E} [\vert \sum_{k \neq i} X_{ik} \vert^2 ]+ \vert z \vert^2\\
& \leq \frac{1}{n} + 1 + \vert z \vert^2 < \infty.
\end{split}
\end{equation*}
We conclude that $\mathbb{E} [\nu_{M -z}]$ is tight. 
\end{proof}

\begin{proof}[Proof of Theorem \ref{thm:singular-value}]
By Lemma \ref{lem:epsilon5}, 
we have that the difference between $a(q)$ and $\mathbb{E} [a (q)] $ converges to zero in probability as $n \rightarrow \infty$. Therefore, in order to prove Theorem \ref{thm:singular-value}, it suffices to 
show the convergence of $\mathbb{E} [a (q)]$ to $\widetilde{\mu}_{|a+g_\gamma-z|}$ for all $q=q(z,\eta) \in \mathbb{H}_+$.

By Theorem \ref{thm:fix-point-equ},  any accumulation point of $\mathbb{E} [\Gamma_n (q)]$ is a solution of the fixed point equation
\begin{equation*}
\begin{bmatrix} \alpha &
\beta \\  \bar{\beta} & \alpha \end{bmatrix} = - \mathbb{E} \left( q+ \begin{bmatrix} \alpha &
\gamma \bar{\beta} \\   \gamma \beta & \alpha \end{bmatrix}- \begin{bmatrix} 0 &
G \\ \overline{G} & 0 \end{bmatrix} \right)^{-1}.
\end{equation*}
By Lemma \ref{lemma:unique}, there is a unique solution of the above fixed point equation for any $q \in \mathbb{H}_+$. In particular, by \eqref{eqn:alpha-Cauchy}, 
\[
   \alpha(q)=s_{A+B, 11}(z,\varepsilon)=\sqrt{-1} \varepsilon \phi\bigg(\big( (z-a-g_\gamma)(z-a-g_\gamma)^*+\varepsilon^2\big)^{-1} \bigg).
\]
By Lemma \ref{lem:tight}, $\alpha(q)$ will be the Stieltjes transform of some symmetric probability measure, denoted by $\widetilde{\nu}_z$. On the other hand, recall that $\widetilde{\mu}$ denotes the symmetrization of a measure $\mu$, hence the right-hand side of the above equation can be rewritten as
\[
 \sqrt{-1} \varepsilon \phi\bigg(\big( (z-a-g_\gamma)(z-a-g_\gamma)^*+\varepsilon^2\big)^{-1} \bigg)=\int_\mathbb{R}\frac{1}{u-\sqrt{-1}\varepsilon}d\widetilde{\mu}_{|a+g_\gamma-z|}(u).
\]
This implies that for any $\varepsilon>0$, two Stieltjes transforms are identical to each other
\[
    S_{\widetilde{\nu}_z}(i\varepsilon)=S_{\widetilde{\mu}_{|a+g_\gamma-z|}}(i\varepsilon). 
\]
We conclude that $\nu_z=\mu_{|a+g_\gamma-z|}$ for any $z\in\mathbb{C}$. 
\end{proof}


\section{Bounds on the least singular value}\label{sec:lsv}
	

	Recall that, for an $n \times n$ matrix $A$, we let $s_1(A) \geq \cdots \geq s_n(A) \geq 0$ be the ordered singular values of $A$.  In this section, we prove the following least singular value bound for $L$. 
		\begin{theorem} \label{thm:lsv}
		Let $(\xi_1, \xi_2)$ be a random vector in $\mathbb{C}^2$ satisfying Assumption \ref{assump:atom}; in addition, assume $\xi_1$ and $\xi_2$ both have unit variance.  Let $X$ be an $n \times n$ random matrix satisfying Conditions {\bf C0-a} and {\bf C0-b} with atom variables $(\xi_1, \xi_2)$, and let $L$ be defined as in \eqref{eqn:defn-Laplacian}.  Fix $z \in \mathbb{C}$ with $z \neq 0$. There exists a constant $c > 0$ (depending only on $z$ and the distribution of $(\xi_1, \xi_2)$) so that 	
\[ \lim_{n \to \infty} \P ( s_n(L-z \sqrt{n}) \leq n^{-c} ) = 0. \]
	\end{theorem}

	We emphasis that Theorem \ref{thm:lsv} does not require the entries to have mean zero.  The rest of this section is devoted to the proof of Theorem \ref{thm:lsv}.

\subsection{Distance of a random vector to a subspace}
		
We start with the following result which follows the arguments from \cite[Lemma 3.5]{BordenaveCaputoChafai2014markov} and \cite[Lemma 7.6]{NguyenORourke2015} with only slight modifications.  
	\begin{lemma}[Distance of a random vector to a subspace] \label{lemma:dist}
	Let $\xi_1$ and $\xi_2$ be complex-valued random variables with unit variance, and let $\psi: \mathbb{N} \to \mathbb{N}$ be such that $\psi(n) < n$ for all $n \geq 2$ and $\psi(n) \to \infty$ as $n \to \infty$.    Then there exists $\eps > 0$ so that the following holds.  
	Let $R = (R_i)_{i=1}^n$ be a random vector with independent coordinates, where for each $1 \leq i \leq n$, $R_i$ is equal in distribution to either $\xi_1$ or $\xi_2$.  
	For any vector $v \in \mathbb{C}^n$ and any subspace $H$ of $\mathbb{C}^n$ with $1 \leq \dim(H) \leq n - \psi(n)$, one has
	\[ \P \left( \dist(R + v, H) \leq \eps \sqrt{n - \dim(H)} \right) \leq C \left( e^{-\eps \psi(n)^3/n^2} + e^{-\eps \psi(n)^2/n} \right), \]
	where $C > 0$ is an absolute constant.  
	\end{lemma}
	\begin{proof}
	By absorbing the expectation $\E[R]$ into the vector $v$, we can assume that $\xi_1$ and $\xi_2$ have mean zero.  In addition, as noted in \cite[Lemma 3.5]{BordenaveCaputoChafai2014markov}, by replacing $H$ with the span of $H$ and $v$ (which only increases the dimension of the subspace by at most one), we can assume that $v = 0$ at the cost of possibly decreasing the value of $\eps$.  Thus, it suffices to prove the lemma under the additional assumptions that $\xi_1$ and $\xi_2$ have mean zero and $v = 0$.  In addition, we assume $\psi(n) \geq C_0 n^{2/3}$ for a constant $C_0 > 0$ to be chosen later as the bound is trivial otherwise.  
	
	We begin with a truncation.  Fix $c' \in (1/2, 1)$.  Let $\E^{(i)}$ denote the conditional expectation with respect to the event $\{ |\xi_i| \leq T\}$ for some parameter $T > 0$ and $i = 1, 2$; note that $\E^{(i)}$ also depends on $T$, but we do not denote this dependence in our notation.  Since $\xi_1$ and $\xi_2$ have unit variance it follows that 
	\[ \min_{i = 1,2} \E^{(i)}[ | \xi_i - \E^{(i)}[\xi_i] |^2 ]  \longrightarrow 1 \]
	as $T \to \infty$.  Thus, we can choose $T_0$ sufficiently large so that
	\begin{equation} \label{eq:choosec}
		\min_{i = 1,2} \E^{(i)}[ | \xi_i - \E^{(i)}[\xi_i] |^2 ] \geq c' \qquad \text{for all } T \geq T_0. 
	\end{equation} 
	Set 
	\begin{equation} \label{eq:trunc}
		T_n := \frac{2n}{c \psi(n)} 
	\end{equation} 
	for some constant $c \in (0, 1/100)$ to be chosen momentarily.  Observe that $T_n > \frac{2}{c}$ for all $n \geq 2$.  Fixing $c$ sufficiently small so that $\frac{2}{c} \geq T_0$, we conclude that
	\begin{equation} \label{eq:T0}
		T_n > \frac{2}{c} \geq T_0 \qquad \text{for all } n \geq 2. 
	\end{equation}
	In order to truncate the entries of $R$, we note that, by Chebyshev's inequality, 
	\begin{equation} \label{eq:vecentrybnd}
		\max_{1 \leq i \leq n} \P (|R_i| > T_n) \leq \frac{1}{T_n^2}. 
	\end{equation}
	Define the event
	\[ E_n = \left\{ \sum_{i=1}^n \un_{\{|R_i| > T_n\}} \geq c \psi(n) \right\}. \]
	By \eqref{eq:vecentrybnd} and Hoeffding's inequality, it follows that
	\[ \P (E_n) \leq \exp\left( - \frac{c^2 \psi(n)^2 }{2n} \right). \]
	Therefore, it suffices to prove the lemma on the event $E_n^c$.  In fact, by conditioning on the indices $i$ for which $|R_i| \leq T_n$, it suffices to prove the result by conditioning on the event 
	\[ \mathcal{E}_m = \{ |R_i| \leq T_n, 1 \leq i \leq m \}, \]
	where $m = n - \lceil c \psi(n) \rceil$.  
	
	On the event $\mathcal{E}_m$, the random vector $R$ may have non-zero mean, so we need to recenter the vector once more.  To this end, let $\E_m$ denote the conditional expectation with respect to the event $\mathcal{E}_m$ and the $\sigma$-algebra $\mathcal{F}_m = \sigma(R_{m+1}, \ldots, R_m)$.  Let $H'$ be the subspace spanned by $H, u, w$, where
	\[ u = (0, \ldots, 0, R_{m+1}, \ldots, R_n), \; w = (\E_m[R_1], \ldots, \E_m[R_m], 0, \ldots, 0). \]
	Clearly, $H'$ is $\mathcal{F}_m$-measurable, and $\dim(H') \leq \dim(H) + 2$. 
	Let $R' = R - u - w$.  It follows that
	\begin{equation} \label{eq:distdist'}
		\dist(R, H) \geq \dist(R, H') = \dist(R', H'). 
	\end{equation} 
	
	It remains to bound $\dist(R', H')$.  By construction, each entry of $R'$ has mean zero.  Since each entry of the original vector $R$ is equal in distribution to either $\xi_1$ or $\xi_2$, it follows from \eqref{eq:choosec} and \eqref{eq:T0} that 
	\begin{equation} \label{eq:sigma}
		\min_{1 \leq i \leq m}  \E_m[ | R'_i|^2 ] \geq c'
	\end{equation} 
	for all $n \geq 2$.  
	By Talagrand's concentration inequality \cite{Talagrand1995concentration}, for any $t > 0$, 
	\begin{equation} \label{eq:talagrand}
		\P_m \left( \left| \dist(R', H') - M_m \right| \geq t \right) \leq 4 \exp \left( - \frac{t^2}{64 T_n^2} \right), 
	\end{equation}
	where $M_m$ is the median of $\dist(R', H')$ under $\mathcal{E}_m$ and $\P_m$ is the conditional probability with respect to $\mathcal{E}_m$ and $\mathcal{F}_m$.  It follows from \eqref{eq:talagrand} (see \cite[Lemma E.3]{TaoVu2010gafa}) that
	\[ M_m \geq \sqrt{ \E_m \dist^2(R', H') } - C T_n, \]
	where $C > 0$ is an absolute constant.  If $P$ denotes the orthogonal projection onto $(H')^\perp$, then
	\begin{align*}
		\E_m \dist^2(R', H') &= \sum_{i=1}^m \E_m [ |R'_i|^2] P_{ii} \\
		&\geq c' \left( \sum_{i=1}^n P_{ii} - \sum_{i=m+1}^n P_{ii} \right) \\
		&\geq c' \left( n - \dim(H') - (n - m) \right)
	\end{align*}
	by \eqref{eq:sigma}.  Taking $0 < \eps < \frac{1}{4} \sqrt{ \frac{49}{50} c'}$, it follows that
	\begin{align*}
		M_m &\geq \sqrt{ c' \left( n - \dim(H') - (n - m) \right)} - CT_n \\
		&\geq  \frac{1}{2}\sqrt{ c' ( n - \dim(H) - 2c \psi(n) )} \\
		&\geq  \frac{1}{2} \sqrt{ \frac{49}{50} c' ( n - \dim(H))} \\
		&\geq 2 \eps \sqrt{n - \dim(H)} 
	\end{align*}
	for $n$ sufficiently large, where we used that $n - \dim(H) \geq \psi(n) \to \infty$, $0 < c < 1/100$, and the fact that $\psi(n) \geq C_0 n^{2/3}$ for a sufficiently large constant $C_0$ (chosen in terms of $c$ and $c'$) to drop the factor of $CT_n$.  Therefore, by \eqref{eq:distdist'} and \eqref{eq:talagrand}, we conclude that
	\begin{align*}
		\P_m \left(\dist(R, H) \leq \eps \sqrt{n - \dim(H)} \right) &\leq \P_m \left( \dist(R', H') \leq M_m - \eps\sqrt{n - \dim(H)} \right) \\
		&\leq \P_m \left( |\dist(R', H') - M_m | \geq \eps \sqrt{\psi(n)} \right) \\
		&\leq 4 \exp \left( - \frac{\eps^2 \psi(n)}{64 T_n^2} \right) \\
		&= 4 \exp \left( - \frac{ \eps' \psi^3(n)}{n^2} \right)
	\end{align*}
	for some constant $\eps' > 0$ depending only on $\eps$ and $c$, and the proof is complete.  
	\end{proof}
	
	For $\delta \in (0,1)$, define the set of sparse vectors 
	\[ \Sparse(\delta) = \{ x \in \mathbb{C}^n : |\supp(x)| \leq \delta n \}, \]
	where $\supp(x) = \{i : x_i \neq 0 \}$ and $|\supp(x)|$ is its cardinality.  Given $\rho \in (0, 1)$ consider the partition of the unit sphere $\mathbb{S}^{n-1}$ into the set of compressible and incompressible vectors:
	\begin{align*}
		\Comp(\delta, \rho) &= \{ x \in \mathbb{S}^{n-1} : \dist(x, \Sparse(\delta)) \leq \rho \} \\
		\Incomp(\delta, \rho) &= \mathbb{S}^{n-1} \setminus \Comp(\delta, \rho). 
	\end{align*}
	For any $n \times n$ matrix $A$, 
	\begin{equation} \label{eq:compincomp}
		s_n(A) = \min_{x \in \mathbb{S}^{n-1}} \|A x \| = \min \left\{ \min_{x \in \Comp(\delta, \rho)} \|Ax\|, \min_{x \in \Incomp(\delta, \rho)} \|A x\| \right \}. 
	\end{equation} 
	
	We will also need the following two results from \cite{RV2008adv} (see also Appendix A in \cite{Bordenave-Chafai-circular}).  
	\begin{lemma} \label{lemma:incomp-cat}
	If $x \in \Incomp(\delta, \rho)$, then there exists a subset $\pi \subset \{1, \ldots, n\}$ so that $|\pi| \geq \delta n/2$ and 
	\[ \frac{\rho}{\sqrt{n}} \leq |x_i| \leq \sqrt{ \frac{2}{\delta n}} \]
	for all $i \in \pi$.  
	\end{lemma}
	
	\begin{lemma} \label{lemma:incomp}
	Let $A$ be any random $n \times n$ matrix, and let $C_k$ denote its $k$-th column.  For $1 \leq k \leq n$, let $H_k = \Span\{ C_j : j \neq k\}$.  Then, for any $t \geq 0$, 
	\[ \P \left( \min_{x \in \Incomp(\delta, \rho)} \|A x\| \leq \frac{t \rho}{\sqrt{n}} \right) \leq \frac{2}{\delta n} \sum_{k=1}^n \P( \dist(C_k, H_k) \leq t ). \]
	\end{lemma}
	
	We will apply \eqref{eq:compincomp} and Lemma \ref{lemma:incomp} to the matrix $A = (L - z \sqrt{n})^{\mathrm{T}}$.  Here, we take the transpose as it will be more convenient to work with rows of $L - z\sqrt{n}$ rather than its columns.

	\subsection{Compressible vectors}
	
	This section is devoted to the following bound.
	\begin{lemma} \label{lemma:comp}
	Let $\kappa > 0$ and $z \in \mathbb{C}$ be fixed.  There exists $\eps > 0$ so that if 
	\[ \delta = \frac{1}{\log n}, \qquad \rho = \frac{\eps}{n^{\kappa} \sqrt{\delta}} \]
	then 
	\[ \lim_{n \to \infty} \P \left( \min_{x \in \Comp(\delta, \rho)} \|(L-z\sqrt{n})^{\mathrm{T}} x \| \leq \frac{\eps}{\sqrt{\delta}}, \|L - z\sqrt{n}\| \leq n^{\kappa} \right) = 0. \]
	\end{lemma}
	\begin{proof}
	If $A$ is an $n \times n$ matrix and $y \in \mathbb{C}^n$ is such that $\supp(y) \subset \pi \subset [n]$, then
	\[ \|A y\| \geq \|y\| s_n(A_{[n] \times \pi}), \]
	where $A_{[n] \times \pi}$ is the $n \times |\pi|$ matrix formed from the columns of $A$ selected by $\pi$.  It follows by the definition of compressible vectors that
	\begin{equation} \label{eq:compbnd}
		\min_{x \in \Comp(\delta \rho)} \|A x \| \geq (1 - \rho) \min_{\pi: |\pi| = \lfloor \delta n \rfloor} s_n(A_{[n] \times \pi}) - \rho \|A\|. 
	\end{equation} 
	To bound the right-hand side, observe that, for $x \in \mathbb{C}^\pi$, 
	\begin{align*}
		\|A_{[n] \times \pi} x \|^2 &= \left\| \sum_{i \in \pi} x_i C_i \right\|^2 \\
		&\geq \max_{i \in \pi} |x_i|^2 \dist^2(C_i, H_i) \\
		&\geq \min_{i \in \pi} \dist^2(C_i, H_i) \frac{1}{|\pi|} \sum_{j \in \pi} |x_j|^2, 
	\end{align*}
	where $C_i$ is the $i$-th column of $A$ and $H_i = \Span\{C_j : j \neq i, j \in \pi \}$.  It follows that
	\begin{equation} \label{eq:distbnd}
		s_n(A_{[n] \times \pi}) \geq \frac{1}{\sqrt{|\pi|}} \min_{i \in \pi} \dist(C_i, H_i). 
	\end{equation} 
	
	We will apply \eqref{eq:distbnd} to $A = (L - z\sqrt{n})^{\mathrm{T}}$.  To this end, fix $\pi \subset [n]$ with $|\pi| = \lfloor \delta n \rfloor$ and fix $i \in \pi$.  
	For any $j \in [n]$, let $R_j$ denote the $j$-th row of $L - z\sqrt{n}$ (i.e., the $j$-th column of $(L-z\sqrt{n})^{\mathrm{T}}$), and let $R'_j$ be the vector $R_j$ with the $i$-th entry removed.  Similarly, define $H_i = \Span \{ R_j : j \neq i, j \in \pi\}$ and $H'_i = \Span \{R'_j : j \neq i, j \in \pi\}$.  It follows that
	\begin{equation} \label{eq:removebnd}
		\dist(R_i, H_i) \geq \dist(R_i', H_i'). 
	\end{equation}
	In particular, by construction, $R_i'$ is an $(n-1)$-vector with independent entries, where each entry is equal in distribution to either $\xi_1$ or $\xi_2$.  In addition, $R_i'$ is independent of $H_i'$.  Thus, by Lemma \ref{lemma:dist} (with $\psi(n) = (n - 1) - 2 \lfloor \delta n \rfloor \geq n/2$), there exists $C, \eps' > 0$ so that
	\[ \P( \dist(R_i', H_i') \leq \eps' \sqrt{n} ) \leq Ce^{-\eps' n}, \]
	and hence
	\[ \P( \dist(R_i, H_i) \leq \eps' \sqrt{n} ) \leq Ce^{-\eps' n} \]
	by \eqref{eq:removebnd}.  Therefore, by \eqref{eq:distbnd} and the union bound, we find
	\[ \P \left( s_n(A_{[n] \times \pi}) \leq \frac{\eps'}{\sqrt{\delta}} \right) \leq C n e^{-\eps' n}. \]
	Returning to \eqref{eq:compbnd}, taking $\eps = \eps'/4$, and applying the union bound (as well as the fact that $1 - \rho \geq 1/2$ for $n$ sufficiently large) gives
	\begin{align} \label{eq:comptozero}
		\P \left( \min_{x \in \Comp(\delta, \rho)} \Vert (L-z\sqrt{n})^{\mathrm{T}}x \Vert \leq \frac{\eps'}{4\sqrt{\delta}}, \|L - z\sqrt{n}\| \leq n^{\kappa} \right) \leq C \binom{n}{\lfloor \delta n \rfloor} n e^{-\eps' n} \to 0
	\end{align}
	as $n \to \infty$.   This completes the proof.  
	\end{proof}
	
	\subsection{Incompressible vectors}
	Let $R_k$ denote the $k$-th row of $L-z\sqrt{n}$ (i.e., $R_k$ is the $k$-th column of $(L - z\sqrt{n})^{\mathrm{T}}$), and let $H_k = \Span \{ R_j : j \neq k\}$.  The goal of this section is to establish the following bound.  
	
	\begin{lemma} \label{lemma:incomp-control}
	Fix $\kappa > 0$ and $z \in \mathbb{C}$ with $z \neq 0$.  
	Let $\delta$ and $\rho$ be as in Lemma \ref{lemma:comp}.  Then there exists constants $C, c > 0$ so that 
	\[ \max_{1 \leq k \leq n} \P \left( \dist(R_k, H_k) \leq c \frac{|z| \rho \delta_0}{n^{\kappa + 1/2}}, \|L - z\sqrt{n}\| \leq n^\kappa  \right) \leq \frac{C}{\sqrt{\eps_0 \delta n }}, \]
	where $\delta_0$ and $\eps_0$ are from Assumption \ref{assump:atom}.  
	\end{lemma}
	
	The proof of Lemma \ref{lemma:incomp-control} is divided into several lemmas.  To start, we fix $k \in [n]$.  Note that all of our bounds below will be uniform in $k$.  Let $\phi$ be the all-ones vector, and set $\Phi = \Span \{ \phi\} = \{ \lambda \phi : \lambda \in \mathbb{C}\}$.  We start with the following deterministic bound.  
	
	\begin{lemma} \label{lemma:normal-dist}
	Let $n \geq 2$.  Any unit vector $\zeta$ orthogonal to $H_k$ satisfies 
	\[ \dist(\zeta, \Phi) \geq \frac{|z|}{\|L - z\sqrt{n}\|}. \]
	In particular, if $\|L - z\sqrt{n}\| \leq s$, one has
	\[ \min_{\lambda \in \mathbb{C}} \| \zeta - \lambda \phi \| \geq \frac{|z|}{s}. \]
	\end{lemma}
	\begin{proof}
	The proof is based on the proof of \cite[Lemma 3.11]{BordenaveCaputoChafai2014markov}.  Set $\hat \phi =  \phi/ \sqrt{n}$ so that $\| \hat \phi \| = 1$.  Since $\hat \phi \in \mathbb{R}^n$, we have
	\[ \dist(\zeta, \Phi) = \dist ( \hat \phi, \Span\{ \zeta\}) = \dist(\hat \phi, \Span\{\bar \zeta\}). \]
	Let $B$ be the matrix obtained from $L - z\sqrt{n}$ by replacing the $k$-th row by zero.  Then $B \bar \zeta = 0$, and so $\bar \zeta \in \ker B$.  It follows that
	\[ \dist(\zeta, \Phi) = \dist(\hat \phi, \Span\{\bar \zeta\}) \geq \dist(\hat \phi, \ker B). \]
	To bound the distance on the right-hand side, observe that we can write $\hat \phi = a u + b v$, where $a, b \in \mathbb{C}$ and $u, v$ are unit vectors with $v \in \ker B$ and $u \in (\ker B)^{\perp}$.  It follows that $|a| = \dist(\hat \phi, \ker B)$.  To bound $|a|$, observe that $(L - z\sqrt{n}) \phi = -z\sqrt{n} \phi$ and so $B \phi = - z\sqrt{n} \phi + z \sqrt{n}e_k$, where $e_k$ is the unit vector with one in the $k$-th coordinate and zeros in all other coordinates.  Thus, we conclude that
	\[|a| \|B\| \geq |a| \|B u \| = \|B \hat \phi \| \geq \sqrt{n-1} |z| \geq |z|. \]
	The conclusion now follows since $\|B \| \leq \|L -z\sqrt{n} \|$.  
	\end{proof}
	
	We will also need the following result.  
	
	\begin{lemma} \label{lemma:normal-not-compressible}
	Fix $\kappa > 0$ and $z \in \mathbb{C}$ with $z \neq 0$.  Let $\delta$ and $\rho$ be as in Lemma \ref{lemma:comp}.  There exists constants $C, \eps > 0$ so that
	\[ \P( \exists \lambda \in \mathbb{C} : \eta(\lambda) \in \Comp(\delta, \rho), \|L - z \sqrt{n} \| \leq n^\kappa) \leq C e^{-\eps n}, \]
	where $\eta(\lambda) = \frac{ \zeta - \lambda \phi }{\| \zeta - \lambda \phi \|}$ and $\zeta$ is any unit vector orthogonal to $H_k$.   
	\end{lemma}
	\begin{proof}
	The proof is based on the proof of \cite[Lemma 3.12]{BordenaveCaputoChafai2014markov}.  
	Note that, by Lemma \ref{lemma:normal-dist}, $\eta(\lambda)$ is well defined since $z \neq 0$, and $\| \eta(\lambda) \| = 1$ for all $\lambda \in \mathbb{C}$.  Let $B$ be the matrix obtained from $L - z\sqrt{n}$ by replacing the $k$-th row by zero.  If there exists $\lambda \in \mathbb{C}$ so that $\eta(\lambda) \in \Comp(\delta, \rho)$, then $B (\bar \eta(\lambda) + \lambda' \phi) = 0$ for $\lambda' = \frac{\bar \lambda}{\| \zeta - \lambda \phi\|} \in \mathbb{C}$.  This implies that
	\[ \min_{x \in \Comp(\delta, \rho), \lambda \in \mathbb{C}} \| B(x + \lambda \phi) \| = 0. \]
	Define $\Phi_k = \Span \{ \phi, e_k \}$, where $e_k$ is the unit vector with one in the $k$-th coordinate and zeros in all other coordinates.  Then $(L - z\sqrt{n}) \phi = -z\sqrt{n} \phi \in \Phi_k$ and $B \phi = -z\sqrt{n} \phi + z\sqrt{n}e_k \in \Phi_k$.  In addition, by construction, $(L - z \sqrt{n})x - Bx \in \Phi_k$ for any vector $x \in \mathbb{C}^n$.  Thus, we conclude that
	\[ \min_{x \in \Comp(\delta, \rho), v \in \Phi_k } \| (L - z\sqrt{n}) x + v \| = 0. \]
	This immediately implies that
	\[ \min_{x \in \Comp(\delta, \rho)} \| \Pi (L - z\sqrt{n}) x \| = 0, \]
	where $\Pi$ is the orthogonal projection onto the orthogonal complement of $\Phi_k$.  
	Therefore, we conclude that
	\begin{align} \label{eq:Picomp}
		\P( \exists \lambda \in \mathbb{C} &: \eta(\lambda) \in \Comp(\delta, \rho), \|L - z \sqrt{n}\| \leq n^\kappa) \\
		&\leq \P\left( \min_{x \in \Comp(\delta, \rho)} \| \Pi (L - z\sqrt{n}) x \| = 0,  \|L - z\sqrt{n} \| \leq n^\kappa\right). \nonumber
	\end{align}
	
	To bound the right-hand side, we apply \eqref{eq:compbnd} and \eqref{eq:distbnd} with $A = \Pi (L - z\sqrt{n})$.  Indeed, we have
	\[ s_n(A_{[n] \times \pi}) \geq \frac{1}{\sqrt{|\pi|}} \min_{i \in \pi} \dist(\Pi C_i, H_i), \]
	where $C_i$ is the $i$-th column of $(L - z\sqrt{n})$ and $H_i = \Span \{ \Pi C_j : j \in \pi, j \neq i\}$.  Fix $\pi \subset [n]$ with $|\pi| = \lfloor \delta n \rfloor$, and fix $i \in \pi$.  Set $H_i' = \Span\{ \Pi C_j, \phi, e_k : j \in \pi, j \neq i\}$.  It follows that
	\[ \dist(\Pi C_i, H_i) \geq \dist(\Pi C_i, H_i') = \dist(C_i, H_i'). \]
	In addition, by construction, $\dim(H_i') \leq \dim(H_i) + 2$.  Unfortunately, $C_i$ and $H_i'$ are dependent.  To get around this, define $I = [n] \setminus (\pi \cup \{i\})$ and let $\Psi: \mathbb{C}^n \to \mathbb{C}^{I}$ be given by $\Psi(v) = v_I = (v_i)_{i \in I}$.  In other words, $\Psi$ removes those entries indexed by $\pi \cup \{i\}$.  Then
	\[ \dist(C_i, H_i') \geq \dist(\Psi(C_i), \Psi(H_i')), \]
	where $\Psi(H_i') = \{\Psi(v) : v \in H_i'\}$.  In addition, $\Psi(C_i)$ is independent of $\Psi(H_i')$.  Thus, by Lemma \ref{lemma:dist} (with $\psi(n) = n - 3 \delta n \geq n/2$), there exists $C, \eps > 0$ so that
	\[ \P \left( s_n(A_{[n] \times \pi} \right) \leq C n e^{-\eps n}. \]
	In particular, as in the proof of \eqref{eq:comptozero}, one finds that
	\[ \P\left( \min_{x \in \Comp(\delta, \rho)} \| \Pi (L - z\sqrt{n}) x \| = 0,  \|L - z\sqrt{n} \| \leq n^\kappa\right) \leq C' \exp(-\eps' n) \]
	for some constants $C, \eps' > 0$.  In view of \eqref{eq:Picomp}, the proof is complete.  
	\end{proof}
	
	We are now in position to complete the proof of Lemma \ref{lemma:incomp-control}.  
	
	\begin{proof}[Proof of Lemma \ref{lemma:incomp-control}]
	Recall that $R_j$ denotes the $j$-th row of $L-z\sqrt{n}$, and let $X_j$ denote the $j$-th row of $X$.  Let $\zeta$ be a unit vector orthogonal to $H_k$.  Clearly, 
	\begin{equation} \label{eq:distnormal}
		\dist(R_k, H_k) \geq | \langle \zeta, R_k \rangle |, 
	\end{equation} 
	where $\langle \cdot, \cdot \rangle$ is the inner product on $\mathbb{C}^n$.  
	By Lemma \ref{lemma:normal-dist}, on the event $\| L - z\sqrt{n} \| \leq n^{\kappa}$, we have
	\begin{equation} \label{eq:normalbnd}
		\min_{\lambda \in \mathbb{C}} \| \zeta - \lambda \phi \| \geq \frac{|z|}{n^{\kappa}}. 
	\end{equation} 
	We now condition on the rows $X_j$, $j \neq k$.  This conditioning fixes $R_j$, $j \neq k$ and $\zeta$.  In fact, by Lemma \ref{lemma:normal-not-compressible}, it suffices to prove the result by conditioning on $X_j$, $j \neq k$ so that \eqref{eq:normalbnd} holds and $\eta(\lambda) \in \Incomp(\delta, \rho)$ for all $\lambda \in \mathbb{C}$ (recall that $\eta(\cdot)$ is defined in the statement of Lemma \ref{lemma:normal-not-compressible}).  While this conditioning fixes $\zeta$, it also affects the entries of $X_k$ due to the dependence between the mirrored entries of $X$.  However, it is straightforward to check that, under this conditioning, the entries of $X_k$ are still independent.  In addition, by condition \eqref{assump:cond}, the entries of $X_k$ are not completely determined by this conditioning.

	Returning to \eqref{eq:distnormal}, we can write
	\[ \langle \zeta, R_k \rangle = \sum_{i=1}^n ( \bar{\zeta}_i - \bar{\zeta}_k ) X_{ki} + \langle \zeta, v \rangle = \langle \zeta - \zeta_k \phi, X_k \rangle + \langle \zeta, v \rangle, \]
	where $v$ is a deterministic vector depending only on $z$ and $X_k$ is the $k$-th row of $X$.  By \eqref{eq:normalbnd},  
	\begin{equation} \label{eq:distalmost}
		\dist(R_k, H_k) \geq \frac{|z|}{n^\kappa} \left| \langle \eta(\zeta_k), X_k \rangle + w \right|, 
	\end{equation}
	where $w$ is depends only on $\zeta$ and $v$ (and hence is fixed by our conditioning).  For notational convenience, let $x = (x_i)_{i \in [n]}$ denote the unit vector $\eta(\zeta_k)$.  Since $x \in \Incomp(\delta, \rho$), Lemma \ref{lemma:incomp-cat} implies the existence of $\pi \subset [n]$ with $|\pi| \geq \delta n/2$ so that 
	\begin{equation} \label{eq:pispread}
		\frac{\rho}{\sqrt{n}} \leq |x_i| \leq \sqrt{ \frac{2}{\delta n}} 
	\end{equation} 
	for all $i \in \pi$.  In view of \eqref{eq:distalmost}, it suffices to upper bound
	\begin{equation} \label{eq:probtobound}
		\sup_{u \in \mathbb{C}}  \widetilde{\mathbb{P}} \left( \left| \sum_{i \in \pi} x_i X_{ki} - u \right| \leq \frac{\delta_0 \rho}{\sqrt{n}} \right), 
	\end{equation} 
	where, by independence, we have conditioned on $X_{ki}$, $i \not\in \pi$ and absorbed their contribution into $u$.  Here, $\widetilde{\mathbb{P}}$ denotes the conditional probability, where we have conditioned on $X_j$, $j \neq k$ and $X_{ki}$, $i \not\in \pi$.  In order to bound \eqref{eq:probtobound}, we apply the Kolmogorov--Rogozin inequality; see \cite{Esseen1968,Rogozin1961}.  In particular, we will apply the multidimensional version \cite[Corollary 1]{Esseen1968}, which implies that
	\[ \sup_{u \in \mathbb{C}} \widetilde{\mathbb{P}} \left( \left| \sum_{i \in \pi} x_i X_{ki} - u \right| \leq \frac{\delta_0 \rho}{\sqrt{n}} \right) \leq C \frac{1}{\sqrt{ \sum_{i \in \pi} (1 - Q_i) }}, \]
	where $C > 0$ is an absolute constant and 
	\[ Q_i = \sup_{u \in \mathbb{C}} \widetilde{\mathbb{P}} \left( |x_i X_{ki} - u| \leq \frac{\delta_0 \rho}{\sqrt{n}} \right). \]
	From \eqref{eq:pispread}, it follows that 
	\[ Q_i \leq \sup_{u \in \mathbb{C}} \widetilde{\mathbb{P}} \left( | X_{ki} - u| \leq \delta_0 \right). \]
	Thus, \eqref{assump:cond} from Assumption \ref{assump:atom} implies that $Q_i \leq 1- \eps_0$ almost surely.  Therefore, we conclude that 
	\[ \sup_{u \in \mathbb{C}} \widetilde{\mathbb{P}} \left( \left| \sum_{i \in \pi} x_i X_{ki} - u \right| \leq \frac{\delta_0 \rho}{\sqrt{n}} \right) \leq \frac{C}{\sqrt{\eps_0 |\pi| }} \]
	almost surely.  Combining this bound with \eqref{eq:distalmost} and the bound $|\pi| \geq \delta n /2$ completes the proof.  
	\end{proof}

	\subsection{Completing the proof of Theorem \ref{thm:lsv}}
	
	We will need the following naive norm bound.  
	\begin{lemma}[Spectral norm bound] \label{lemma:norm}
	Fix $z \in \mathbb{C}$.  For any $p > 0$, there exists $C, \kappa > 0$, so that 
	\[ \P ( \| L - z\sqrt{n} \| > n^{\kappa}) \leq C n^{-p}. \]
	\end{lemma}
	\begin{proof}
	Recall that $L = X - D$.  Since $D$ contains the row sums of $X$, it suffices to assume that the diagonal entries of $X$ are zero.  As $z$ is fixed, it follows that for $\kappa > 1/2$
	\begin{equation} \label{eq:sufficesnorm}
		\P( \| L - z\sqrt{n} \| > n^{\kappa}) \leq \P ( \| X \| > n^{\kappa - 1}) + \P( \| D \| > n^{\kappa -1 }) 
	\end{equation}
	for all $n$ sufficiently large.  Thus, it suffices to bound the two probabilities on the right-hand side.  Since $\xi_1$ and $\xi_2$ have finite second moment, there exists $C_0 > 0$ so that
	\begin{equation} \label{eq:chebyentry}
		\max_{i=1, 2} \P( |\xi_i| > t ) \leq \frac{C_0}{t^2} 
	\end{equation} 
	for all $t > 0$.  Define the event
	\[ \Omega = \left\{ \max_{i,j \in [n]} |X_{ij}| \leq n^{1 + p/2} \right\}. \]
	It follows from \eqref{eq:chebyentry} and the union bound that 
	\begin{equation} \label{eq:Omegac}
		\P(\Omega^c) \leq \frac{ C_0 n^2 }{n^{2 + p}} = \frac{C_0}{n^{p}}. 
	\end{equation} 
	On the event $\Omega$, 
	\[ \|X\|^2 \leq \|X\|^2_2 = \sum_{i,j \in [n]} |X_{ij}|^2 \leq n^{4 + p}, \]
	where $\|X\|_2 = \sqrt{\Tr (X X^\ast)}$  is the Frobenius norm of $X$.  In addition, since $D$ contains the row sums of $X$, by the triangle inequality, 
	\[ \|D\| = \max_{i} |D_{ii}| \leq n^{2 + p/2} \]
	on the event $\Omega$.  Thus, in view of \eqref{eq:sufficesnorm} and \eqref{eq:Omegac}, taking $\kappa = 3 + p/2$ completes the proof.  
	\end{proof}
	
	We are now ready to prove Theorem \ref{thm:lsv}.  
	
	\begin{proof}[Proof of Theorem \ref{thm:lsv}]
	Fix $z \in \mathbb{C}$ with $z \neq 0$.  Let $\delta$ and $\rho$ be as in Lemma \ref{lemma:comp}.  From Lemma \ref{lemma:norm}, there exists $C, \kappa > 0$ so that
	\begin{equation} \label{eq:norm100}
		\P ( \| L - z\sqrt{n} \| > n^{\kappa}) \leq C n^{-100}. 
	\end{equation} 
	Combining \eqref{eq:norm100} with Lemma \ref{lemma:comp} shows there exists $\eps > 0$ so that 
	\begin{equation} \label{eq:comp-conc}
		\lim_{n \to \infty} \P \left( \min_{x \in \Comp(\delta, \rho)} \|(L-z\sqrt{n})^{\mathrm{T}} x \| \leq \frac{\eps}{\sqrt{\delta}} \right) = 0. 
	\end{equation}
	Combining \eqref{eq:norm100} with Lemma \ref{lemma:incomp-control} implies the existence of $C, c > 0$ so that
	\[ \max_{1 \leq k \leq n} \P \left( \dist(R_k, H_k) \leq c \frac{|z| \rho \delta_0}{n^{\kappa + 1/2}} \right) \leq \frac{C}{\sqrt{\eps_0 \delta n }}, \]
	where $\delta_0$ and $\eps_0$ are from Assumption \ref{assump:atom}.  Thus, in view of Lemma \ref{lemma:incomp}, we conclude that
	\begin{align}
		\P \left( \min_{x \in \Incomp(\delta, \rho)} \|(L -z\sqrt{n})^{\mathrm{T}} x\| \leq c \frac{|z| \rho^2 \delta_0}{n^{\kappa + 1}}  \right) &\leq \frac{2}{\delta n} \sum_{k=1}^n \P \left( \dist(R_k, H_k) \leq c \frac{|z| \rho \delta_0}{n^{\kappa + 1/2}}  \right) \nonumber \\
		&\leq \frac{2C}{\delta \sqrt{\eps_0 \delta n }} \longrightarrow 0 \label{eq:incomp-conc}
	\end{align}
	as $n \to \infty$.  Thus, combining \eqref{eq:comp-conc} and \eqref{eq:incomp-conc} with \eqref{eq:compincomp} completes the proof of the theorem.  
	\end{proof}

	\section{Moderately small singular values}\label{sec:mslv}
	
	In this section, we prove the following lower bound on the moderately small singular values of $L$.  
	\begin{lemma}\label{lem:moderate}
	Let $X$ be as in Theorem \ref{thm:lsv}, and fix $z \in \mathbb{C}$.  Let $a_n = n/\log^2 n$.  Then there exists $\eps > 0$ so that almost surely, for $n$ sufficiently large, 
	\[ s_{n - i}\left( M - z \right) \geq \eps \frac{i}{n} \]
	for all $a_n \leq i \leq n-1$.  
	\end{lemma}
	\begin{proof}
	The proof is based on the proofs of \cite[Lemma 3.14]{BordenaveCaputoChafai2014markov} and \cite[Lemma 7.6]{NguyenORourke2015}, which in turn are based on the arguments of Tao and Vu \cite{TaoVu2010}.  Let $s_1 \geq \cdots \geq s_n$ denote the singular values of $M - z$, and fix $a_n \leq i \leq n - 1$.  Let $B$ be the matrix formed from the first $m = n - \lfloor i/2 \rfloor$ rows of $L - z \sqrt{n}$.  Let $s_1' \geq \cdots \geq s_m'$ denote the singular values of $B$.  By Cauchy's interlacing property, it follows that
	\[ \frac{1}{\sqrt{n}} s_{n - i}' \leq s_{n - i}. \]
	By \cite[Lemma A.4]{TaoVu2010}, 
	\[ s_1'^{-2} + \cdots + s_{m}'^{-2} = \dist_1^{-2} + \cdots + \dist_m^{-2}, \]
	where $\dist_j = \dist(R_j, H_j)$ is the distance from the $j$-th row $R_j$ of $B$ to $H_j$, the subspace spanned by all the other rows of $B$, i.e., $H_j = \Span\{ R_i : 1 \leq i \leq m, i \neq j\}$.  It follows that
	\begin{equation} \label{eq:singvaluedist}
		\frac{i}{2n} s_{n-i}^{-2} \leq \frac{i}{2} s_{n-i}'^{-2} \leq \sum_{j = n - i}^m s_{j}^{-2} \leq \sum_{j=1}^m \dist_{j}^{-2}. 
	\end{equation} 
	We now estimate $\dist_j$, $1 \leq j \leq m$.  Unfortunately, $R_j$ and $H_j$ are dependent.  To work around this problem, fix $1 \leq j \leq m$, and define the matrix $B_j$ to be the matrix $B$ with the $j$-th column removed.  Let $R_i'$ be the $i$-th row of $B_j$, and set $H_j' = \Span \{ R_i' : 1 \leq i \leq m, i \neq j\}$.  Since 
	\[ \dist_j \geq \dist(R_j', H_j'), \]  
	it suffices to bound $\dist(R_j', H_j')$.  Notice that the entries of $R_j'$ are independent and each entry is equal in distribution to either $\xi_1$ or $\xi_2$.  In addition, $R_j'$ is independent of $H_j'$, so we can apply Lemma \ref{lemma:dist}.  In this case, $\dim(H_j') \leq n - i/2 \leq n - 1 - \frac{1}{4} a_n$.  Applying Lemma \ref{lemma:dist} (with $\psi(n) = \frac{1}{4} a_n$), we find constants $C, \eps > 0$ so that
	\[ \sup_{1 \leq j \leq m} \P \left( \dist_j \leq \eps \sqrt{i} \right) \leq C \exp \left( - \eps \frac{n}{\log^6 n} \right). \]
	By the union bound, we have 
	\[ \sum_{n=1}^\infty \P \left( \bigcup_{a_n \leq i \leq n-1} \bigcup_{j = 1}^m \left\{ \dist_j \leq \eps \sqrt{i} \right\} \right) < \infty, \]
	and hence, by the Borel--Cantelli lemma, we find that almost surely, for $n$ sufficiently large, 
	\[ \dist_j \geq \eps \sqrt{i} \]
	for all $1 \leq j \leq m$ and all $a_n \leq i \leq n-1$.  The proof of the lemma is now complete by combining the bound above with \eqref{eq:singvaluedist}.  
	\end{proof}

\section{Convergence of the empirical spectral distribution}
In this section, we will prove Theorem \ref{thm:eigenvalue}.
We will use the following standard Hermitization method \cite[Lemma 4.3]{Bordenave-Chafai-circular}:
\begin{lemma}[Hermitization]\label{lem:Hermitian}
Let $(X_n)_{n \geq1}$ be a sequence of $n \times n$ complex random matrixes. Suppose that there exists a family of probability measures $(\nu_z)_{z \in \mathbb{C}}$ on $\mathbb{R}_+$ such that for a.a. $z \in \mathbb{C},$
\begin{enumerate}[\rm (i)]
\item $\nu_{X_n -z}$ tends weakly in probability to $\nu_z;$
\item $\log (\cdot)$ is uniformly integrable in probability for $(\nu_{X_n -z})_{n \geq 1}.$
\end{enumerate} 
Then, in probability, $\mu_{X_n}$ converges weakly to the probability measure $\mu$ given by
\begin{equation*}
\mu= \frac{1}{2\pi} \triangle \int_0^\infty \log (t) d\nu_z(t).
\end{equation*} 
\end{lemma}

\begin{theorem}[Uniform Integrability] \label{thm:uniform-inter}
Under the assumptions of Theorem \ref{thm:eigenvalue}, there exists an increasing function $J: [0, \infty) \rightarrow [0, \infty)$ with $J(t)/t \rightarrow \infty$ as $t \rightarrow \infty$
such that for a.a. $z \in \mathbb{C},$
\begin{equation}
\lim_{t \rightarrow \infty} \limsup_{n \rightarrow \infty} \P \left( \int_0^\infty J (\vert \log (s) \vert) d \nu_{M-z} (s) >t \right) =0.
\end{equation} 
\end{theorem}

\begin{proof}
We follow the proof of \cite[Theorem 3.2]{BordenaveCaputoChafai2014markov}. Choose $J(t) = t^2.$ By Lemma \ref{lem:tight}, we have 
\begin{equation*}
\mathbb{E} \int s^2 d\nu_{M-z}(s)< \infty. 
\end{equation*}
By the Markov inequality, there exists a constant $C>0$ such that for any $t>0,$
\begin{equation*}
\P \left( \int s^2 d\nu_{M-z}(s) > t \right) \leq C t^{-1}.
\end{equation*}
Hence, it suffices to show that 
\begin{equation}\label{eq:truncation}
\lim_{t \rightarrow \infty} \limsup_{n \rightarrow \infty} \P \left( \int_0^1 J (\vert \log (s) \vert) d \nu_{M-z} (s) >t \right) =0.
\end{equation} 

Denote $n_\ast$ by the largest $i$ such that $s_{n-i}(M-z) \leq 1.$ With the notation of Lemma \ref{lem:moderate}, let $F_n$ be the event that $s_n(M-z) \geq n^{-c -1/2}$ and that for all
$a_n \leq i \leq n-1, s_{n-i} (M-z) \geq \varepsilon i/n.$ Then by Theorem \ref{thm:lsv} and Lemma \ref{lem:moderate}, $F_n$ holds with probability tending to $1$ as $n \rightarrow \infty.$ Therefore, on the event $F_n$, one has 
\begin{equation*}
\begin{split}
\frac{1}{n} \sum_{i=1}^{n_\ast} J(\vert \log s_{n-i} \vert ) & \leq \frac{1}{n} \sum_{i=1}^{a_n} J (\vert \log (n^{-c -1/2}) \vert ) + \frac{1}{n} \sum_{i=1}^n J (\vert \log (\varepsilon i/n)\vert)\\
& = \frac{ \log^2(n^{-c -1/2})}{\log^2 n} + \frac{1}{n} \sum_{i=1}^n  \log^2(\varepsilon i/n) = O(1), 
\end{split}
\end{equation*}
where the sum is uniformly bounded as it can be identified with a Riemann sum of a finite integral.  
Indeed, we have proved that there exists a constant $C>0$ such that
\begin{equation}
\lim_{n \rightarrow \infty} \P \left( \frac{1}{n} \sum_{i=1}^{n_\ast} J (\vert \log s_{n-i}(M-z) \vert) \leq C \right)=1,
\end{equation}
which concludes the proof of \eqref{eq:truncation}. 
\end{proof}

\begin{proof}[Proof of Theorem \ref{thm:eigenvalue}] 
We apply Lemma \ref{lem:Hermitian} to our matrix $M.$ The condition {\rm (i)} follows from Theorem \ref{thm:singular-value}. On the other hand, from the de la Vall\'{e}e-Poussin criterion for uniform integrability, Theorem \ref{thm:uniform-inter} implies the condition {\rm (ii)}, i.e., the uniform integrability in probability of $\log (\cdot)$ for measures $(\nu_{M-z})_{n \geq 1}$.
Hence, the limit distribution $\mu$ is given as follows:
\begin{equation*}
\mu = \frac{1}{2\pi} \triangle \int \log (t) d\nu_z(t).
\end{equation*}

Note that for any $x \in \log^+ (\mathcal{A}),$
\begin{equation*}
\phi( \log|x- z|) = \int \log (t) d\mu_{|x-z|} (t).
\end{equation*}
Hence, the Brown measure of $x$ is given by the formula
\begin{equation*}
\mu_x = \frac{1}{2\pi} \triangle \int \log (t) d\mu_{|x-z|} (t).
\end{equation*}
Recall that by \eqref{eqn:nu-1.12},
\[
    \nu_z=\mu_{|a+g_\gamma-z|},\qquad \forall z\in\mathbb{C}.
\]
Therefore, the limit distribution $\mu$ is the same as the Brown measure of $a+g_\gamma$. 
\end{proof}

\bigskip

{\bf Acknowledgment.} S. O'Rourke was supported in part by the National Science Foundation under Grant No. DMS-2143142.  Z.Yin was partially supported by NSFC No. 12031004. P. Zhong was supported in part by a grant from Simons Foundation and NSF award LEAPS-MPS-2316836.

\bibliographystyle{acm}
\bibliography{Elliptic-law}

\begin{thebibliography}{10}

\bibitem{Gernotbook}
{\sc Akemann, G., Baik, J., and Di~Francesco, P.}
\newblock {\em {The Oxford Handbook of Random Matrix Theory}}.
\newblock Oxford University Press, 09 2015.

\bibitem{MR4388923}
{\sc Alt, J., and Kr\"{u}ger, T.}
\newblock Local elliptic law.
\newblock {\em Bernoulli 28}, 2 (2022), 886--909.

\bibitem{book2010RMT}
{\sc Anderson, G.~W., Guionnet, A., and Zeitouni, O.}
\newblock {\em An introduction to random matrices}, vol.~118 of {\em Cambridge
  Studies in Advanced Mathematics}.
\newblock Cambridge University Press, Cambridge, 2010.

\bibitem{BaiBook}
{\sc Bai, Z., and Silverstein, J.~W.}
\newblock {\em Spectral Analysis of Large Dimensional Random Matrices}.
\newblock Springer Series in Statistics. Springer New York, NY, 2006.

\bibitem{Baron2023}
{\sc Baron, J.~W.}
\newblock A path integral approach to sparse non-{H}ermitian random matrices.
\newblock {\em arXiv:2308.13605\/} (2023).

\bibitem{BR2019}
{\sc Basak, A., and Rudelson, M.}
\newblock {The circular law for sparse non-Hermitian matrices}.
\newblock {\em The Annals of Probability 47}, 4 (2019), 2359 -- 2416.

\bibitem{BelinschiTR2018-sub-operator-valued}
{\sc Belinschi, S.~T., Mai, T., and Speicher, R.}
\newblock Analytic subordination theory of operator-valued free additive
  convolution and the solution of a general random matrix problem.
\newblock {\em J. Reine Angew. Math. 732\/} (2017), 21--53.

\bibitem{BelinschiYinZhong2021Brown}
{\sc Belinschi, S.~T., Yin, Z., and Zhong, P.}
\newblock Brown measure of the sum of a triangular elliptic operator and a free
  random variable.
\newblock {\em arXiv:2209.11823v2\/} (2022).

\bibitem{MR3744883}
{\sc Benaych-Georges, F., C\'{e}bron, G., and Rochet, J.}
\newblock Fluctuation of matrix entries and application to outliers of elliptic
  matrices.
\newblock {\em Canad. J. Math. 70}, 1 (2018), 3--25.

\bibitem{BordenaveCaputoChafai2014markov}
{\sc Bordenave, C., Caputo, P., and Chafa\"{\i}, D.}
\newblock Spectrum of {M}arkov generators on sparse random graphs.
\newblock {\em Comm. Pure Appl. Math. 67}, 4 (2014), 621--669.

\bibitem{Bordenave-Chafai-circular}
{\sc Bordenave, C., and Chafa\"{\i}, D.}
\newblock Around the circular law.
\newblock {\em Probab. Surv. 9\/} (2012), 1--89.

\bibitem{BrycDemboJiang2006}
{\sc Bryc, W.~o., Dembo, A., and Jiang, T.}
\newblock Spectral measure of large random {H}ankel, {M}arkov and {T}oeplitz
  matrices.
\newblock {\em Ann. Probab. 34}, 1 (2006), 1--38.

\bibitem{MR4565612}
{\sc Byun, S.-S., Kang, N.-G., Lee, J.~O., and Lee, J.}
\newblock Real eigenvalues of elliptic random matrices.
\newblock {\em Int. Math. Res. Not. IMRN}, 3 (2023), 2243--2280.

\bibitem{MR4489824}
{\sc Campbell, A., and O'Rourke, S.}
\newblock Spectrum of heavy-tailed elliptic random matrices.
\newblock {\em Electron. J. Probab. 27\/} (2022), Paper No. 125, 56.

\bibitem{MR1484954}
{\sc Chang, J.~T., and Pollard, D.}
\newblock Conditioning as disintegration.
\newblock {\em Statist. Neerlandica 51}, 3 (1997), 287--317.

\bibitem{MR3930614}
{\sc Durrett, R.}
\newblock {\em Probability---theory and examples}, fifth~ed., vol.~49 of {\em
  Cambridge Series in Statistical and Probabilistic Mathematics}.
\newblock Cambridge University Press, Cambridge, 2019.

\bibitem{Esseen1968}
{\sc Esseen, C.~G.}
\newblock On the concentration function of a sum of independent random
  variables.
\newblock {\em Z. Wahrscheinlichkeitstheorie und Verw. Gebiete 9\/} (1968),
  290--308.

\bibitem{MR4201596}
{\sc Fyodorov, Y.~V., and Tarnowski, W.}
\newblock Condition numbers for real eigenvalues in the real elliptic
  {G}aussian ensemble.
\newblock {\em Ann. Henri Poincar\'{e} 22}, 1 (2021), 309--330.

\bibitem{MR2999219}
{\sc Girko, V.}
\newblock The elliptic law. {T}hirty years later.
\newblock {\em Random Oper. Stoch. Equ. 20}, 4 (2012), 347--399.

\bibitem{MR3068415}
{\sc Girko, V.}
\newblock The generalized elliptic law.
\newblock {\em Random Oper. Stoch. Equ. 21}, 2 (2013), 191--215.

\bibitem{MR0781018}
{\sc Girko, V.~L.}
\newblock An elliptic law.
\newblock {\em Dokl. Akad. Nauk Ukrain. SSR Ser. A}, 1 (1985), 56--59.

\bibitem{MR0816278}
{\sc Girko, V.~L.}
\newblock The elliptic law.
\newblock {\em Teor. Veroyatnost. i Primenen. 30}, 4 (1985), 640--651.

\bibitem{MR0901506}
{\sc Girko, V.~L.}
\newblock Elliptic law and elements of {$G$}-analysis.
\newblock In {\em Probability theory and mathematical statistics, {V}ol. {I}
  ({V}ilnius, 1985)}. VNU Sci. Press, Utrecht, 1987, pp.~489--507.

\bibitem{MR1354817}
{\sc Girko, V.~L.}
\newblock The elliptic law: ten years later. {I}.
\newblock {\em Random Oper. Stochastic Equations 3}, 3 (1995), 257--302.

\bibitem{MR1373142}
{\sc Girko, V.~L.}
\newblock The elliptic law: ten years later. {II}.
\newblock {\em Random Oper. Stochastic Equations 3}, 4 (1995), 377--398.

\bibitem{MR1483014}
{\sc Girko, V.~L.}
\newblock Strong elliptic law.
\newblock {\em Random Oper. Stochastic Equations 5}, 3 (1997), 269--306.

\bibitem{MR2215672}
{\sc Girko, V.~L.}
\newblock The strong elliptic law. {T}wenty years later. {I}.
\newblock {\em Random Oper. Stochastic Equations 14}, 1 (2006), 59--102.

\bibitem{MR3403838}
{\sc G\"{o}tze, F., Naumov, A., and Tikhomirov, A.}
\newblock On a generalization of the elliptic law for random matrices.
\newblock {\em Acta Phys. Polon. B 46}, 9 (2015), 1737--1745.

\bibitem{HWineq}
{\sc Hoffman, A.~J., and Wielandt, H.~W.}
\newblock {The variation of the spectrum of a normal matrix}.
\newblock {\em Duke Math. J. 20}, 1 (1953), 37 -- 39.

\bibitem{MNR-JPA}
{\sc Metz1, F.~L., Neri, I., and Rogers, T.}
\newblock Spectral theory of sparse non-{H}ermitian random matrices.
\newblock {\em J. Phys. A: Math. Theor. 52\/} (2019), 434003.

\bibitem{MingoSpeicherBook}
{\sc Mingo, J.~A., and Speicher, R.}
\newblock {\em Free probability and random matrices}, vol.~35 of {\em Fields
  Institute Monographs}.
\newblock Springer, New York; Fields Institute for Research in Mathematical
  Sciences, Toronto, ON, 2017.

\bibitem{Naumov13}
{\sc Naumov, A.}
\newblock Elliptic law for real random matrices.
\newblock {\em arXiv:1201.1639\/} (2012).

\bibitem{NM-PRL}
{\sc Neri, I., and Metz, F.~L.}
\newblock Spectra of sparse non-{H}ermitian random matrices: An analytical
  solution.
\newblock {\em Phys. Rev. Lett. 109\/} (Jul 2012), 030602.

\bibitem{NguyenORourke2015}
{\sc Nguyen, H.~H., and O'Rourke, S.}
\newblock The elliptic law.
\newblock {\em Int. Math. Res. Not. IMRN}, 17 (2015), 7620--7689.

\bibitem{2013Low}
{\sc O'Rourke, S., and Renfrew, D.}
\newblock Low rank perturbations of large elliptic random matrices.
\newblock {\em Electron J Probab 19}, 43 (2014), 1--65.

\bibitem{MR3540493}
{\sc O'Rourke, S., and Renfrew, D.}
\newblock Central limit theorem for linear eigenvalue statistics of elliptic
  random matrices.
\newblock {\em J. Theoret. Probab. 29}, 3 (2016), 1121--1191.

\bibitem{MR3357969}
{\sc O'Rourke, S., Renfrew, D., Soshnikov, A., and Vu, V.}
\newblock Products of independent elliptic random matrices.
\newblock {\em J. Stat. Phys. 160}, 1 (2015), 89--119.

\bibitem{Resnick}
{\sc Resnick, S.~I.}
\newblock {\em A Probability Path}.
\newblock Modern Birkh\"{a}user Classics. Birkh\"{a}user Boston, MA, 2013.

\bibitem{Rogozin1961}
{\sc Rogozin, B.~A.}
\newblock On the increase of dispersion of sums of independent random
  variables.
\newblock {\em Teor. Verojatnost. i Primenen 6\/} (1961), 106--108.

\bibitem{RotmanBook}
{\sc Rotman, J.~J.}
\newblock {\em An Introduction to Algebraic Topology}.
\newblock Graduate Texts in Mathematics. Springer New York, NY, 1988.

\bibitem{RV2008adv}
{\sc Rudelson, M., and Vershynin, R.}
\newblock The {L}ittlewood-{O}fford problem and invertibility of random
  matrices.
\newblock {\em Adv. Math. 218}, 2 (2008), 600--633.

\bibitem{SpeicherAMS1998}
{\sc Speicher, R.}
\newblock Combinatorial theory of the free product with amalgamation and
  operator-valued free probability theory.
\newblock {\em Mem. Amer. Math. Soc. 132}, 627 (1998), x+88.

\bibitem{Talagrand1995concentration}
{\sc Talagrand, M.}
\newblock Concentration of measure and isoperimetric inequalities in product
  spaces.
\newblock {\em Inst. Hautes \'{E}tudes Sci. Publ. Math.}, 81 (1995), 73--205.

\bibitem{TaoVu2010gafa}
{\sc Tao, T., and Vu, V.}
\newblock Random matrices: the distribution of the smallest singular values.
\newblock {\em Geom. Funct. Anal. 20}, 1 (2010), 260--297.

\bibitem{TaoVu2010}
{\sc Tao, T., and Vu, V.}
\newblock Random matrices: universality of {ESD}s and the circular law.
\newblock {\em Ann. Probab. 38}, 5 (2010), 2023--2065.
\newblock With an appendix by Manjunath Krishnapur.

\bibitem{Tikhom2023}
{\sc Tikhomirov, A.~N.}
\newblock Limit theorem for spectra of {L}aplace matrix of random graphs.
\newblock {\em Mathematics 11}, 3 (2023).

\bibitem{Voiculescu1991}
{\sc Voiculescu, D.}
\newblock Limit laws for random matrices and free products.
\newblock {\em Invent. Math. 104}, 1 (1991), 201--220.

\bibitem{DVV-operator-valued-1992}
{\sc Voiculescu, D.}
\newblock Operations on certain non-commutative operator-valued random
  variables.
\newblock {\em Ast\'{e}risque}, 232 (1995), 243--275.
\newblock Recent advances in operator algebras (Orl\'{e}ans, 1992).

\bibitem{Zhong2021Brown_ctgt}
{\sc Zhong, P.}
\newblock Brown measure of the sum of an elliptic operator and a free random
  variable in a finite von neumann algebra.
\newblock {\em arXiv:2108.09844v4\/} (2021).

\end{thebibliography}

\end{document}